\documentclass[10pt]{article}
\usepackage[latin2]{inputenc}
\usepackage{amsmath}
\usepackage[pdftex]{graphicx}
\usepackage{amsfonts}
\usepackage{amssymb}
\usepackage{algorithm}
\usepackage{algorithmic}
\usepackage[extra]{tipa}
\usepackage[left=1.5in]{geometry}
\usepackage{geometry}
\usepackage{wrapfig}
\usepackage{setspace}
\usepackage{epsfig}
\usepackage{color}
\usepackage{longtable}
\usepackage{listings}
\usepackage{tikz}
\usepackage{cancel}
\usepackage{float}
\usepackage{cite}
\usepackage{url}
\usepackage{hyperref}
\usepackage{amssymb, amsmath, epsfig, euscript, graphicx, xspace, amsthm, amscd, amsmath, mathtools}
\usepackage{float}
\restylefloat{table}
\usepackage{amsthm}
\usepackage{array,float}
\usepackage{setspace}
\usepackage{calligra}
\usepackage{enumitem}
\usepackage{tikz}
\usetikzlibrary{decorations.pathreplacing,calligraphy}
\usepackage{graphicx}
\usepackage{fancyref}
\usepackage{caption}
\usetikzlibrary{arrows,automata,positioning}
\usetikzlibrary{arrows}

\begin{document}

\newcommand{\Div}{\text{Div}}
\newcommand{\Pic}{\text{Pic}}
\newcommand{\Span}{\text{Span}}
\newcommand{\Edg}{\text{Edg}}
\newcommand{\ZZ}{\mathbb{Z}}
\newcommand{\QQ}{\mathbb{Q}}
\newcommand{\RR}{\mathbb{R}}
\newcommand{\CC}{\mathbb{C}}
\newcommand{\FF}{\mathbb{F}}
\newcommand{\NN}{\mathbb{N}}
\newcommand{\Cal}{\mathcal}
\newcommand{\mubar}{\overline{\mu}}
\newcommand{\bye}{\end{document}}
\newcommand{\aenum}{\begin{enumerate}[label=(\alph*)]}
\newcommand{\nenum}{\begin{enumerate}[label=(\arabic*)]}
\newcommand{\renum}{\begin{enumerate}[label=(\roman*)]}
\def\nor{\trianglelefteq}
\def\iso{\cong}
\def\Aut{\text{Aut}}
\def\QQ{\Bbb Q}
\def\ZZ{\Bbb Z}
\def\RR{\Bbb R}
\def\CC{\Bbb C}
\def\NN{\Bbb N}
\def\FF{\Bbb F}
\def\ord#1{| \, #1 \, |}
\def\gp#1{\langle \, #1 \, \rangle}
\def\Div{\text{Div}}
\def\Pic{\text{Pic}}
\def\vPic{\text{vPic}}
\def\Pr{\text{Pr}}
\def\RB{R\mathcal B}
\def\Gal{\text{Gal}}
\def\nor{\trianglelefteq}
\def\iso{\cong}
\def\Aut{\text{Aut}}
\def\QQ{\Bbb Q}
\def\ZZ{\Bbb Z}
\def\RR{\Bbb R}
\def\CC{\Bbb C}
\def\NN{\Bbb N}
\def\FF{\Bbb F}
\def\ord#1{| \, #1 \, |}
\def\gp#1{\langle \, #1 \, \rangle}
\def\Div{\text{Div}}
\def\Pic{\text{Pic}}
\def\vPic{\text{vPic}}
\def\Pr{\text{Pr}}
\def\RB{R\mathcal B}
\def\Gal{\text{Gal}}
\def\Char{\text{Char}}
\def\bigmid{\, \big | \,}
\renewcommand{\epsilon}{\varepsilon}
\newcommand{\DD}{\mathbb{D}}
\newtheorem{Theorem}{Theorem}
\newtheorem{Lemma}{Lemma}
\newtheorem{Corollary}{Corollary}
\newtheorem{proposition}{Proposition}
\newtheorem{definition}{Definition}
\theoremstyle{definition}
\newtheorem{exmp}{Example}

\begin{center}
\textbf{THE JACOBIAN OF CYCLIC VOLTAGE COVERS OF $K_n$}\\
\vspace{10px}
\textbf{SOPHIA R GONET}\\
\end{center}

\noindent \textbf{Abstract}
 This paper proves results about the Jacobians of a certain family of covering graphs, $Y$, of a base graph $X$, that is constructed from an assignment of elements from a group $G$ to the edges of $X$ ($G$ is called the {\it voltage group} and $Y$ is called the {\it derived graph}\/). Of particular interest is when the voltage assignment is given by mapping a generator of the cyclic group of order $d$ to a single edge of $X$ (all other edges are assigned the identity), called a {\it single voltage assignment}. 
%The voltage group $G$ acts as graph automorphisms of the derived graph $Y$, the group of divisors of $Y$ becomes a module over the group ring $\ZZ[G]$, and the Laplacian endomorphism on the group of divisors of $Y$---which is used to compute the Jacobian of $Y$---can be described by a matrix with entries from $\ZZ[G]$, called the {\it voltage Laplacian}.  Using this and matrix computations
Both the order and abelian group structure of the Jacobian of single voltage assignment derived graphs are determined when the base graph $X$ is the complete graph on $n$ vertices, for every $n$ and $d$. 
%When $G$ is abelian, the determinant of the \textit{voltage Laplacian} matrix is called the {\it reduced Stickelberger element}; and it is shown to be a power of two times the graph Stickelberger element defined in the literature in terms of Ihara zeta-functions.  
Using zeta-functions, general product formulas that relate the order of the Jacobian of $Y$ to that of $X$ are developed; these formulas become very simple and explicit in the special case of single voltage covers of $X$.

\section{Introduction}\label{intro}

%The Jacobian is a finite abelian group whose decomposition into invariant factor cyclic subgroups is its {\it Smith Normal Form}; and the product of these invariant factors, which is the determinant of the reduced Laplacian matrix, is the order of the Jacobian of $X$, the important {\it tree number} of the graph $X$.
The exact structure of the Jacobian is known for only a few classes of graphs. One such family of graphs for which the Jacobian is known is the complete graph on $n$ vertices, $K_n$ (see \cite{Hop14}).  For other papers containing calculations of Jacobians for families of graphs see \cite{BAI2003251, Berget2009TheCG, CHRISTIANSON2002233, Jacobson2003231, DUCEY2014316, Chandler2014, Sin2018227, DUCEY2018154, DUCEY2021105424}. The purpose of this paper is to start with any connected base graph $X$ \footnote[1]{In this paper the term graph will mean a simple graph with no loops or multiple edges.} and construct certain $d$-fold covering graphs $Y$ of $X$ such that the Jacobian of $Y$, denoted by $\mathcal{J}(Y)$, is expressible in terms of $d$ and $\mathcal{J}(X)$, i.e., we obtain formulas for the order of $\mathcal{J}(Y)$; but for some $X$, we also obtain the abelian group structure of $\mathcal{J}(Y)$.  In this way, starting with a graph whose Jacobian is known, we exhibit an infinite family of new graphs whose Jacobians are explicitly and computationally easily determined.

More specifically, we focus on \textit{derived graphs} $Y$ obtained from {\it voltage assignments} to the edges of a base graph $X$ (as described in Section~2)  especially ones in which the voltage assignment is the \textit{single voltage assignment} (given by mapping a generator of the cyclic group of order $d$ to a single edge of $X$ with all other edges assigned the identity). In Section~3, even more specifically, we establish the structure and the order of the Jacobian of single voltage covers of the complete graph on $n$ vertices: \begin{Theorem}\label{KnJac}
Let $Y$ be a single voltage cover of the complete graph $K_n$ by the cyclic group of order $d$, 
where $n \ge 4$ and $d \ge 3$.  Then 
$$
J(Y) \iso (\ZZ/n\ZZ)^{(n-4)d+2} \oplus (\ZZ/n(n-2)\ZZ)^{d-2}\oplus \ZZ/dn(n-2)\ZZ,
$$
where the exponents indicate the multiplicities of the (distinct) invariant factors. In particular, the order of the Jacobian of $Y$ is $n^{(n-3)d+1}\cdot (n-2)^{d-1}\cdot d$ and its rank is $(n - 3)d + 1$.
\end{Theorem}
We also explicitly describe the small cases when $d=2$ and when $n=3$.

In Section \ref{4} we develop some general methods for relating the order of $\mathcal J(Y)$ to $\mathcal J(X)$, using {\it Ihara zeta-functions}. This results in the following Theorem:
\begin{Theorem}\label{OrdJac}
Let $X$ be any connected graph with the number of its edges is not equal to the number of its vertices, and let $Y$ be the derived graph resulting from any voltage assignment of elements of a finite abelian group $G$ to edges of $X$ such that $Y$ is connected. Then the order of the Jacobian of the derived graph $Y$ is
$$
|\mathcal{J}(Y)|=\frac{1}{d}\cdot|\mathcal{J}(X)|\prod_{\chi_t\neq \chi_0}\chi_t(\Theta_{Y/X}).
$$
where $\Theta_{Y/X}$ is a a certain explicitly defined element of the integral group ring $\ZZ[G]$ --- called the reduced Stickelberger element --- and the $\chi_t$ run over all non-principal irreducible characters of $G$.
\end{Theorem}
These formulas result in explicit (closed form) formulas 
for general (connected) base graphs $X$ and arbitrary (connected) single voltage derived graphs $Y$ in terms of the 
reduced Stickelberger element, namely:
\begin{Corollary}\label{OrdJac2}
Let $X$ and $Y$ satisfy the hypotheses of Theorem 1, but let $G$ be the cyclic group of order $d$ and assume the voltage assignment is the single voltage assignment. Let $K$ be the greatest common divisor of the (integer) coefficients of $\Theta_{Y/X}$. Then
$$
|\mathcal{J}(Y)|=|\mathcal{J}(X)|\cdot |K|^{d-1}d.
$$
\end{Corollary}
In the case of single voltage assignments,  we shall see that $\Theta_{Y/X}$ has an especially simple form and $K$ is very easy to compute.

This paper gives a computationally effective method for computing the Smith Normal Form of the Laplacian of the derived graph $Y$. We first do elementary row and column operations on the \textit{voltage} Laplacian, and then we tensor the entries with the regular representation of the group $G$ to obtain the ordinary Laplacian for the derived graph $Y$. We further reduce by row and column operations over $\ZZ$ to obtain the Smith Normal Form for the Laplacian of $Y$. In the midst of manipulations of the voltage Laplacian, after reducing it to nearly diagonal
(or nearly upper-triangular) form, we compute the reduced Stickelberger element (the determinant of the voltage Laplacian when $G$ is abelian).

The terminology and theory of Jacobians, voltage graphs and their derived covering graphs is first summarized in Section \ref{2}. We prove Theorem \ref{KnJac} in Section \ref{3}. In addition, we obtain partial results for the Jacobian of cyclic voltage covers of $K_{n,n}$ in Section \ref{3.5}. In particular, we obtain the primes $p$ (and their powers) that divide the order of the Jacobian of the derived graph $Y$, for $p$ not dividing $n$.

This paper comprises Chapter 3 and 4 of the author's dissertation \cite{phdthesis}, which contains significantly more examples and an array of additional theoretical and computational material on voltage graphs and their associated derived graphs.

\section{Preliminaries}\label{2}
The \textit{Picard} and \textit{Jacobian} groups of a graph $X$ are invariants that may be defined in terms of the \textit{Laplacian} matrix of $X$.  In particular, the Picard group is the cokernel of the Laplacian $L$, considered as a $\ZZ$-module endomorphism of the $\ZZ$-module of divisors of $X$, $\Div_{\ZZ}(X)$, which is free over $\ZZ$ of rank $n = |V(X)|$. Then for the \textit{voltage graph} with \textit{voltage group} $G$ and \textit{derived graph} $Y$, the $\ZZ$-module of divisors of $Y$ becomes a free module of the same rank $n$ over the larger ring $\ZZ[G]$.  The voltage adjacency matrix captures the ``voltage adjacencies'' across the different sheets of the covering of $X$ by $Y$; and these sheets are acted on by $G$ as the regular representation.  The \textit{voltage Laplacian} $\mathcal L$ is then an endomorphism of the free $\ZZ[G]$-module $\Div_{\ZZ}(Y) = \Div_{\ZZ[G]}(X)$, which agrees with the ordinary Laplacian endomorphism for $Y$; it can be represented by an $n \times n$ matrix with entries in $\ZZ[G]$, or an $n\ord G \times n\ord G$ matrix of integers---the ordinary Laplacian matrix for $Y$---by forming $\mathcal L \otimes \rho$, where $\rho$ is a matrix representation for the regular representation of $G$. 

\subsection{The Jacobian group}\label{2.1}
See \cite{corry_perkinson_2018} for further details on this section.

\begin{definition}
A \textit{divisor} on a graph $X$ (possibly infinite) is an element of the free abelian group on the vertices $V=V(X)$:
$$
\Div(X) = \{ \sum_{v \in V(X)}a_v v \mid a_v \in \ZZ \}
$$
where each $\sum_{v \in V(X)} a_v v$ is a formal linear combination of the vertices of $X$ with
integer coefficients, where only finitely many $a_v$ are nonzero (in the case when $X$ is an infinite graph).
When $V(X) = \{ v_1,\dots, v_n\}$, we may write the elements of $\Div(X)$ as 
$a_1 v_1 + a_2 v_2 + \cdots + a_n v_n$, where each $a_i \in \ZZ$.
\end{definition}

\begin{definition}
Let $G=(V,E)$ be a graph with vertices $\{v_1,\cdots, v_n\}.$ The graph \textit{Laplacian} $L=L_X$ is the $n\times n$ matrix given by
$$
L_{i,j}=\begin{cases}
\deg(v_i) &\text{ if }i=j\\
-1 &\text{ if } \text{ $v_i$ is adjacent to $v_j$}\\
0& \text{ if }\text{ $i\neq j$ and $v_i$ is not adjacent to $v_j$}
\end{cases}
$$
\end{definition}

\noindent The Laplacian is also the matrix representation of the following group homomorphism $\mathcal{L}$ defined below. For each fixed $v_i$ define the \textit{principal divisor} $p_i$ based at $v_i$ by
$$
p_i=\deg(v_i)v_i-\sum_{\substack{j=1\\ j\neq i}}^n \delta_{i,j}v_j
$$
where $\delta_{i,j}=1$ if $v_j$ is adjacent to $v_i$ and 0 otherwise. Then
$$
\mathcal{L}:\Div(X)\to \Div(X)\qquad \text{ where }\qquad \mathcal{L}(v_i)=p_i.
$$
When extended by $\ZZ$-linearity to all of $\Div(X)$, this is a $\ZZ$-linear homomorphism from $\Div(X)$ to itself, whose image is $\Pr(X)$, the group of \textit{principal divisors}. From this, we get the definition of the \textit{Picard group}:
$$
\Pic(X)=\Div(X)/\text{im}(\mathcal{L})=\text{coker}(\mathcal L).
$$
(For further details on this, see \cite{corry_perkinson_2018}, Section 2.1.)
%Evidently, the matrix representation for $\mathcal{L}$ with respect to the basis of vertices is $L,$ the Laplacian matrix.\\
A \textit{reduced Laplacian} $\tilde{L}$ is the $(n-1)\times (n-1)$ integer matrix obtained by removing the row and column corresponding to any vertex $v$ from the Laplacian matrix $L$. So the \textit{Jacobian group} can be computed as the cokernel of the reduced Laplacian matrix
$$
\mathcal{J}(X)\cong \ZZ^{n-1}/\text{im}(\widetilde{L})=\text{coker}(\widetilde{L}),
$$
where $\ZZ^{n-1}$ denotes the free $\ZZ$-module on the set $V(X)-\{v\}$ of rank $n-1.$	

\subsubsection{Computing the Jacobian}\label{2.3}

We give a computationally effective method for computing the Jacobian of a graph. We first begin by stating the Structure Theorem for Finitely Generated Abelian Groups, which can be found in \cite{corry_perkinson_2018} as Proposition 2.23. 

\begin{Theorem}\label{structurefg}[Structure Theorem for Finitely Generated Abelian Groups]{}
A group is a finitely generated abelian group if and only if it is isomorphic to
\begin{equation}\label{fg}
\ZZ/{d_1}\ZZ \oplus \cdots \oplus \ZZ/{d_k}\ZZ \oplus \ZZ^r
\end{equation}
for some unique integers $d_1,\cdots, d_k$ with $d_i>1$ for all $i$ and some integer $r\geq 0$ that satisfy the following condition: $d_i\mid d_{i+1} \ \forall i.$ The $d_i$ are the invariant factors of the group and $r$ is called the free rank of the group.
\end{Theorem}		

\noindent If a matrix $M$ may be obtained from a matrix $N$ through a sequence of \textit{elementary integer row and column operations}, we write $M\sim N.$ The next proposition can be found in \cite{corry_perkinson_2018} as Proposition 2.28.

\begin{proposition}
Let $M$ and $N$ be $m\times n$ integer matrices. If $M\sim N$, then $\text{coker}(M)\cong\text{coker}(N).$
\end{proposition}

\begin{definition}
An $m\times n$ integer matrix $M$ is in \textit{Smith Normal Form} (SNF) if\\ $M~=\text{diag}(d_1,\cdots d_k,0,\cdots, 0)$, a diagonal matrix, where $d_1,\cdots, d_k$ are positive integers such that $d_i\mid d_{i+1} \ \forall i.$ The $d_i$ with $d_i\geq2$ are called the invariant factors of $M.$
\end{definition}

\noindent We now illustrate the relevance of the Structure Theorem and Smith Normal Form in computing the Jacobian of graph. The next proposition can be found in \cite{klivans_2019} as Proposition 4.5.2.

\begin{proposition}
If $M$ is a non-singular $n\times n$ integer matrix and the Smith normal form of $M$ is $\text{diag}(d_1,\cdots, d_k)$ then
$$
\text{coker}(M)\cong\bigoplus_{i=1}^k \ZZ/d_i\ZZ.
$$
\end{proposition}

\noindent So for any finite \textit{connected} graph $X$ with reduced Laplacian $\widetilde{L}$ (with respect to any vertex), since $\mathcal{J}(X)=\text{coker}(\widetilde{L}),$ it follows that the invariant factors of $\mathcal{J}(X)$ are determined by the invariant factors of $\widetilde{L}.$ Since $\widetilde{L}$ is invertible, none of its invariant factors are 0. So the free rank of $\mathcal{J}(X)$ is 0. The SNF factors of $\widetilde{L}$ equal to 1 have no effect on the isomorphism class of $\text{coker}(\widetilde{L}).$ Suppose $d_1\cdots, d_k$ are the invariant factors of $\widetilde{L}$ (the SNF factors that are greater than 1). The $d_i$ are the same as the invariant factors of the finite abelian group $\mathcal{J}(X)$ as in Theorem \ref{structurefg}, and hence,
$$
\mathcal{J}(X)\cong \ZZ/{d_1}\ZZ \oplus \cdots \oplus \ZZ/{d_k}\ZZ.
$$
The structure of $\Pic(X)$ is then determined since $\Pic(X)\cong \ZZ\oplus \mathcal{J}(X)$ (by Proposition 1 of \cite{corry_perkinson_2018}): the free rank of $\Pic(X)$ is 1 and 
$$
\Pic(X)\cong \ZZ\oplus \ZZ/{d_1}\ZZ\oplus\cdots \oplus \ZZ/{d_k}\ZZ.
$$
The next proposition can be found in \cite{corry_perkinson_2018} as Proposition 2.37.
\begin{proposition}
For $X$ a connected graph, the order of the Jacobian of $X$ is
$$
|\mathcal{J}(X)|=\det(\widetilde{L}).
$$
for any (hence every) reduced Laplacian $\widetilde{L}$ of $X$.
\end{proposition}

\subsection{Voltage Graphs}\label{2.4}

We first give the definition of a general covering graph.
%We begin this section by defining a covering graph in general. 
%In Section \ref{3.1} we then go on to describe a specific type of covering graph, called a \textit{derived graph}, that arises from what is called a \textit{voltage graph}---where elements from a group (which may be finite or infinite) are assigned to the edges of a fixed base graph $X$.  In Section \ref{3.4}, we give the definition of an \textit{intermediate covering graph} $\widetilde{X}$. We then state several theorems that are important for Section \ref{sec 4}.  In Section \ref{3.5}, for $Y$ a derived graph (which may be infinite), we consider \textit{the group of divisors of $Y$}, \textit{the group of principal divisors of $Y$}, \textit{the Picard group of $Y$}, and \textit{the Laplacian of $Y$}, all as $\ZZ$ and $\ZZ[G]$ modules, where $G$ is any finite or infinite group.  \\

\begin{definition}
An undirected graph $Y$ is a \textit{covering} of an undirected graph $X$ if, after arbitrarily directing the edges of $X$, there is an assignment of directions to the edges of $Y$ and an onto graph homomorphism $\pi:Y\to X$ sending neighborhoods of $Y$ one-to-one onto neighborhoods of $X$ which preserve directions. We call such $\pi$ a covering map.
\end{definition}

\begin{definition}\label{dsheeted}
A \textit{$d$-sheeted} covering means every fiber contains exactly $d$ elements, i.e.,
 $$
 |\pi^{-1}(x)|=d \ \forall x\in V(X).
 $$
\end{definition}

	%\noindent We first start with an undirected, finite and connected graph $X$ and a group $G$ (which may be finite or infinite). We assume, for convenience that $X$ has no loops and no multiple edges; but the discussion easily generalizes to multigraphs. Next, we arbitrarily orient the edges of $X$. Then we label the forward-directed edges of $X$ with elements from the group $G$. (Note: the edges of $X$ do not have to labeled with different group elements. For instance, each edge could be assigned the group identity element.) These labels are referred to as the \textit{voltages} and the assignment itself is called the \textit{voltage assignment}. Furthermore, if we have a directed edge which goes from $u$ to $v$ labeled with the group element $\tau$, then label the edge which goes from $v$ to $u$ by the group element's inverse, namely $\tau^{-1}:$
\begin{definition} Let $X$ be a graph whose edges have been oriented, and let $G$ be a group (finite or infinite). For a fixed orientation of the edges of $X$, let $E(X)^+$ denote the set of forward-directed edges of $X$; and let $E(X)^-$ denote the same edges but each with the reverse orientation (so each undirected edge of $X$ becomes two edges in the disjoint union of $E(X)^+$ and $E(X)^-$). An (ordinary) \textit{voltage assignment} is a map
$$
\alpha: E(X)^+\cup E(X)^-\to G
$$
such that if $e_{i,j}\in E(X)^+$ and $\alpha(e_{i,j})=\alpha_{i,j}\in G,$ then $e_{j,i}\in E(X)^-$ and $\alpha(e_{j,i})=\alpha_{i,j}^{-1}$ (the inverse group element), where $e_{i,j}$ denotes the directed edge from $v_i$ to $v_j.$ 
The triple $(X, G, \alpha)$ is called an (ordinary) \textit {voltage graph}. The values of $\alpha$ are called the \textit{voltages} and $G$ is called the \textit{voltage group}.
Note that a voltage assignment $\alpha$ is uniquely determined by its values on $E(X)^+$, so we will henceforth only specify $\alpha$ on the forward-directed edges of $X$.
\end{definition}

The vertices of $X$ are labeled as $v_1,\dots,v_n$.
This imposes a natural lexicographic orientation on $X$, namely whenever there is an edge between $v_i$ and $v_j$,
orient the edge $v_i \rightarrow v_j$ if $i < j$ (called the \textit{standard orientation}). Note that results on derived graphs do not depend on the choice of orientation by \cite{phdthesis}, so without further mention, we adopt the standard orientation.

Any such voltage assignment can be codified by its $n \times n$ {\it voltage adjacency matrix}: $A_\alpha = \begin{pmatrix} \alpha_{i,j} \end{pmatrix}$
where the $i,j$ entry is zero if there is no edge between $v_i$ and $v_j$; and the diagonal entries (for our graphs) are zero. (Note that the {voltage adjacency matrix} is also defined in \cite{chepuri2019arborescences} as Definition 2.14.)
%The entries of $A_\alpha$ are from the integral group ring $\ZZ [G]$, and $A_\alpha$
%is ``inverse symmetric'' in the sense that its transpose, $A_\alpha^t$, is the same as $A_\alpha$, but with all group elements in nonzero entries inverted. 
%(As usual, we identify the integer $k$ with the group ring element $k \tau_0$, where, for formality, $\tau_0$ is the identity of $G$.)
%Changing the orientation of $X$ simply results in a voltage assignment where some $\alpha_{i,j}$ are inverted. \\

The purpose of assigning voltages to the graph $X$, called \textit{the base graph}, is to obtain an object called \textit{the derived graph}, called $Y$ here. To get the vertices of $Y$, make $d=|G|$ copies of each vertex $x\in V(X)$ labeling them as $x_{\tau_0}, x_{\tau_1}, x_{\tau_2},\cdots,x_{\tau_{d-1}} $where $\tau_0,\tau_1, \tau_2,\cdots, \tau_{d-1} \in G$ (again, the same formal construction works even if $|G|$ is uncountable). So there are $|G|\cdot |V(X)|$ vertices in $Y$. 
Now create the edges of $Y$ by the following rule:
whenever there is an edge from $v_i$ to $v_j$ in the base graph $X$ with assigned voltage $\alpha_{i,j}$, 
create edges that go from $v_{i,g}$ to $v_{j, g \alpha_{i,j}}$ in $Y$, for every $g \in G$,
where $g \alpha_{i,j}$ is the group-product of these two group elements in $G$. 
If $|G|=d,$ then $\pi:Y\to X$ is a $d$-sheeted covering map (where again, $d$ may be any infinite cardinal too). Note that the degree (valence) of each vertex $v_\tau$ of $Y$ is the same as the degree of $v=\pi(v_\tau)$ in $X$. Also, no two vertices in the same fiber of $\pi$ are adjacent in~$Y$.\\

\subsubsection{Constructing the Laplacian for a Derived Graph $Y$}
\label{2.4.1}
This subsection describes a computationally effective way of constructing the adjacency matrix, and hence, Laplacian for the derived graph $Y$
from a given voltage graph $(X,G,\alpha)$. 
%This construction is described in the language of $\ZZ[G]$-modules in Section \ref{2.7}. 
Note that $G$ may be any finite group, not necessarily abelian. 

Since group-voltages ``multiply on the right'' on vertices of $Y$, we view the Laplacian operator as acting on the right on divisors, $\vec{a}L,$ and so we compute the Laplacian matrix from this perspective. Fix a listing of the vertices of $X$ as $v_1,\cdots, v_n$ and let $G= \{\tau_0,\tau_1,\dots,\tau_{d-1}\}.$ 
The right regular representation of $G$ is the homomorphism $\rho : G \rightarrow GL_d(\QQ)$, where for each $g \in G$ the $d \times d$
matrix $\rho(g)$ has a 1 in position $i,j$ if and only if $ \tau_i g = \tau_j$; all other entries are 0. Note that if $G$ is abelian then the left and right regular representations are the same, and so $\rho$ is just called the regular representation.

Fix any ordering of the group elements as $\tau_0,\cdots, \tau_{d-1}.$ Next, for each $i$ list the vertices of $Y$ in the fiber over $v_i$ as
$v_{i,\tau_0}, v_{i,\tau_1},\cdots, v_{i,\tau_{d-1}}.$
Finally, list all the vertices of $Y$ (using the \textit{standard orientation}) with $v_{i,\tau_j}$ before $v_{p,\tau_q}$ if $i<p$ or if $i=p$ with $j<q.$

With respect to this ordering, the adjacency matrix of $Y$ is the ``tensor product'' of $A_\alpha$ and $\rho$ as follows:
Create the $nd \times nd$ (block) matrix by replacing each nonzero entry $\alpha_{i,j}$ in $A_\alpha$ by the $d \times d$ matrix $\rho(\alpha_{i,j})$;
replace each zero entry in $A_\alpha$ by the $d \times d$ zero matrix. 
Denote this matrix by $A_{Y}$. 
%The fact that $A_Y$ is an adjacency matrix for $Y$ is the observation that, 
%by definition of the derived graph, the adjacency matrix for $Y$ can be written as 
%an $n \times n$ matrix whose entries are $\ord G \times \ord G$ block matrices. 
%When $v_i$ is not adjacent to $v_j$, the $i,j$ block is identically zero.
%When $v_i$ is adjacent to $v_j$, i.e., there is a 1 in the $i,j$ entry of $A_X$, then 
%$v_{i,\tau}$ is adjacent to $v_{j,\tau \alpha_{i,j}}$, for every $\tau \in G$;
With respect to the above labeling of the elements of $G$, (which we chose to be the same for all blocks)
the $i,j$ block of an adjacency matrix for $Y$ is a matrix for $\rho(\alpha_{i,j})$ (again, keeping in mind that the matrix acts on the right on row vectors). This method is discussed briefly in \cite{Mizuno1995295}, as well as in \cite{feng_kwak_lee_2004}.

Since each vertex $v_\sigma$ in $Y$ has the same degree as $v$ in $X$, to create the $nd\times nd$ degree matrix of $Y$, likewise replace each entry $n_{i,j}$ of the degree matrix of $X$ by the scalar matrix $n_{i,j}I_d$,
where $I_d$ is the $d \times d$ identity matrix (noting that all off diagonal entries are thus replaced by the zero matrix).
Denote this (diagonal) matrix by~$D_{Y}$. 
The Laplacian of $Y$ is then $D_{Y} - A_{Y}$; 
and the reduced Laplacian, $\widetilde{L}_{Y}$, is obtained from it by deleting the $i^{\text{th}}$ row and column for any $i$.

\subsection{The Voltage Laplacian and the Reduced Stickelberger Element }\label{2.5}

We next clarify the relationship between the $\ZZ$-module of all divisors on $Y$ (where $G$, hence also $Y$ may be infinite) and its structure as a $\ZZ[G]$-module.
%In each of the subsections below we describe the ``$\ZZ$ structure'' as item (1) and the ``$\ZZ[G]$ structure'' of the same object(s) as item (2). 
From this, we define \textit{the reduced Stickelberger element}, and describe how it relates to the Stickelberger element defined in \cite{HMSV19}.

\medskip\noindent

\noindent {\it The group of divisors of $Y$:}
By definition, the divisor group of $Y$, denoted as $\text{Div}(Y)$, is the free $\ZZ$-module on the vertices of $Y$.
Since the vertices are $\{ v_{i,\sigma} \mid 1 \le i \le n, \; \sigma \in G \}$, these form a $\ZZ$-basis of $\ZZ$-rank $n|G|$.
Now let $g\in G,$ $v_i\in V(X)$ and let $\sigma\in G$ be a sheet index. Then
$$
g: v_{i,\sigma}\mapsto v_{i,g\sigma}
$$
defines a left group action on the set of vertices in each fiber of $\pi$. These disjoint orbits are
$\pi^{-1}(v_i) = \{ v_{i,\tau} \ | \ \tau\in G \}$ for $1 \le i \le n$. 
Since $G$ acts as the left regular representation on each fiber, the $\ZZ$-span of each one is a free $\ZZ [G]$-module of rank~1; 
and $\Div(Y)$ is then a direct sum of these: a free $\ZZ[G]$-module of rank~$n$, i.e.
\begin{equation}\label{equation 2.4}
\Div(Y)= \ZZ[G] v_{1,\tau_0}\oplus \ZZ[G]v_{2,\tau_0} \oplus \cdots \oplus \ZZ[G] v_{n,\tau_0}
\end{equation}
where we, for simplicity,  choose the identity representative from each orbit.\\
\\
\noindent {\it The group of principal divisors of $Y$:} To consider principal divisors in $Y$ we need to find what vertices are adjacent to a given $v_{i,{\tau}}$ in $Y$. We have that $v_{i,\tau}$ is adjacent to $v_{j,\delta}$ in $Y$ if and only if $v_i$ is adjacent to $v_j$ in $X$ and one of the following holds: $v_i\xrightarrow{\alpha_{i,j}} v_j$ with $\tau\alpha_{i,j}=\delta,$ or $v_j\xrightarrow{\alpha_{j,i}} v_i$ with $\delta\alpha_{j,i}=\tau.$ Since $\alpha_{j,i}^{-1}=\alpha_{i,j}$ in $G$, in both cases $\delta=\tau\alpha_{i,j}.$ 
This says (independent of the orientation on edges): 
the unoriented degree of each vertex $v_{i,\tau}$ in $Y$ is the same as the unoriented degree of $v_i$ in $X$ (call it $n_i$), and
the principal divisor ``based at $v_{i,\tau}$'' is, by definition,
$$
p_{i,\tau} = n_i \; v_{i,\tau} - \sum_{\substack{j=1 \\ v_i\sim v_j}}^n v_{j,\tau\alpha_{i,j}}.
$$
Then, by definition, the $\ZZ$-module of principal divisors, denoted by $\Pr(Y)$ is
$$
\Pr(Y) = \Span_{\ZZ} \{ p_{i,{\tau}} \mid 1 \le i \le n, \; \tau\in G \}.
$$
Note that these principal generators are not necessarily a $\ZZ$-{\it basis} of $\Pr(Y)$.
Any graph automorphism permutes principal divisors, so $\Pr(Y)$ is a $\ZZ [G]$-module (not generally a free module).
The left action of $G$ partitions the above set of principal divisors into orbits 
$\mathcal O_i = \{ p_{i,\tau} \mid \tau\in G\}$, for $i=1,2,\dots,n$,
and again $G$ acts as the regular permutation representation on each $\mathcal O_i$ 
(it acts by left multiplication on the subscript $\tau$ for $\tau\in G$).
Thus if we pick a representative of each orbit, say for convenience $p_{i,{\tau_0}}$, where $\tau_0$ is the identity of $G$, then we get that 
$$
\Pr(Y) = \Span_{\ZZ [G]} \{  p_{i,\tau_0} \mid 1 \le i \le n \}.
$$

\noindent {\it The Picard group of $Y$:} Since the Picard group is $\Pic(Y) = \Div(Y) / \Pr(Y)$, all of these terms are both $\ZZ$ and $\ZZ [G]$-modules, as described above.\\
\\
\noindent {\it The Laplacian of $Y$:} By definition, if $A_{Y}$ is the (ordinary) adjacency matrix for $Y$ and $D_{Y}$ is the degree matrix, 
then the Laplacian is $L_{Y} = D_{Y} - A_{Y}$, which, when $|G|=d$ is finite, is an $nd \times nd$ matrix with entries from $\ZZ$.
By definition, $L_Y$ can also be written as a $\ZZ$-module endomorphism of $\Div_\ZZ(Y)$ (even when $|G|=\infty$):
$$
L_{Y}(v_{i,\tau}) = p_{i,\tau}, \qquad \text{for all } 1\leq i \leq n \text{ and } \tau \in G
$$ 
and this is extended by $\ZZ$-linearity to all of $\Div(Y)$.
In particular, the $\ZZ$-module image of $\Div(Y)$ under the $\ZZ$-module homomorphism $L_{Y}$ is $\Pr(Y)$, 
and its cokernel is $\Pic(Y)=\Div(Y) / \Pr(Y)$. 
Consider the $\ZZ [G]$-module homomorphism which is defined on the above $\ZZ [G]$-basis in (\ref{equation 2.4}) of $\Div(Y)$ by
$$
\mathcal L_\alpha : \Div(Y) \longrightarrow \Div(Y)\qquad\text{by}\qquad
\mathcal L_\alpha(v_{i,\tau_0}) = p_{i,\tau_0}, \quad 1 \le i \le n.
$$
This map is now extended by $\ZZ [G]$-linearity to all of $\Div(Y)$, namely for all $\tau \in G$ by
$$
\mathcal L_\alpha ( \tau \cdot v_{i,\tau_0}) = \tau \cdot \mathcal L_\alpha (v_{i,\tau_0}) = \tau \cdot p_{i,\tau_0} =p_{i,\tau\cdot  \tau_0}= p_{i,\tau}
$$
and likewise for sums and differences of these.
\begin{definition}
The $\ZZ[G]$-module homomorphism $\mathcal L_\alpha$ is called the \textit{voltage Laplacian} of $Y$.
\end{definition}
\noindent Thus by the action of $G$ we have
$$
\mathcal L_\alpha ( v_{i,\tau}) = p_{i,\tau} \qquad \text{for all } \tau \in G.
$$
Since $G$ acts transitively on the $i^{\text{th}}$ $G$-orbit of both the vertices and the principal divisors in $\Div(Y)$ we see that
$$
\text{\it the image of the $\ZZ [G]$-module homomorphism $\mathcal L_\alpha$ is the $\ZZ [G]$-submodule $\Pr(Y)$.}
$$
Finally, we compute the $n \times n$ matrix of $\mathcal L_\alpha$ with respect to the above $\ZZ [G]$-basis of $\Div(Y)$.
In the $j^{\text{th}}$ column of a matrix and $i^{\text{th}}$ row we put the coefficient of $v_{i,\tau_0}$ in the expansion of $\mathcal L_\alpha(v_{j,\tau_0})$.  By the above, this $i,j$-entry equals $n_i$ if $i=j$ and $-\alpha_{j,i}$ if $i\neq j$
(and zero if $v_i$ is not adjacent to $v_j$ in $X$).
Now this results in the $i,j$-entry of the \textit{transpose} of $D_X-A_\alpha.$ However, since the voltage ``multiplies'' on the right on divisors, we should be computing a matrix $M$ representing $\mathcal{L}_\alpha$ where $\vec{a}M=\vec{b}$ (i.e., a row vector $\times M = $ a row vector). Such $M$ is the transpose of the matrix we just computed (i.e., the matrix we computed acts on the left on column vectors). 
%\noindent This discussion shows that the voltage Laplacian endomorphism of $\Div(Y)$ is exactly the same as the 
%ordinary Laplacian when we consider $\Div(Y)$ as a $\ZZ$-module rather than a $\ZZ[G]$-module.
%But note that the voltage Laplacian is always represented by an $n \times n$ matrix (with entries from $\ZZ[G]$), 
%even when $G$ is an infinite group; whereas the ordinary Laplacian of $Y$ has a matrix representation 
%(of degree $n \ord G$) only when $G$ is finite, as described explicitly in Section 2.2.
We summarize this discussion by the following theorem.

\bigskip\noindent
\begin{Theorem}\label{Stick1}
Let $G$ be any group. The voltage Laplacian $\mathcal L_\alpha : \Div(Y) \longrightarrow \Div(Y)$ is a $\ZZ [G]$-module homomorphism whose image is $\Pr(Y)$ and 
cokernel is $\Pic(Y)$, and its $n \times n$ matrix with respect to the $\ZZ [G]$-basis $v_{1,\tau_0},\dots,v_{n,\tau_0}$ is equal to 
$D_X - A_\alpha$, where $\tau_0$ is the identity of $G$, $D_X$ is the degree matrix for the base graph $X$ and $A_\alpha$ is the voltage adjacency matrix of $X$. 
\end{Theorem}
\begin{definition} 
Assume $G$ is abelian. We call $\Theta_{Y/X}=\det(D_X-A_\alpha)$ the \textit{reduced Stickelberger element.}
\end{definition}
\noindent Note that $\Theta_{Y/X}$ is an element of the integral group ring $\ZZ[G].$ We need $G$ to be commutative only for the determinant to be well-defined (over a commutative ring).\\
\\
\noindent In \cite{HMSV19}, their Stickelberger element is denoted  by $\theta^*_{Y/X}(1) e$.
By their Theorem 4.5, since their map $\phi$ is seen to be the same as our map $\mathcal L_\alpha$,
their Stickelberger element relates to ours by
$$
\theta^*_{Y/X}(1) e = 2^{r_X - 1} \Theta_{Y/X} ,
$$
so our version of the Stickelberger element eliminates a 2-power factor (with nonnegative exponent).  This leads to a stronger annihilation statement than their Theorem 4.7 (and explains the terminology ``reduced''). 
%To be more precise, a careful reading of the {\it proof} of their annihilation result 
%(Theorem~4.7), which relies only on {\it exactly the same} (adjoint matrix) result 
%that we use, reveals that they too actually achieve {\it the same} as the following 
%corollary; however it is {\it stated} in that paper with the power of 2 still in place.
%Again, we use the expression ``reduced Stickelberger element'' merely to 
%distinguish it from the version in \cite{HMSV19}.
\begin{Corollary}\label{annih}
For $G$ abelian, $\Theta_{Y/X}$ annihilates the Picard group $\Pic(Y)$
of any derived graph (viewed as a $\ZZ$-module or a $\ZZ [G]$-module), hence it also annihilates $\mathcal{J}(Y).$
\end{Corollary}
\begin{proof}
This is immediate from Exercise 3 in Section 11.4 of \cite{DF04} applied to $\mathcal L_\alpha$, viewing $\Pic(Y)$ as a 
module over the commutative ring $\ZZ[ G]$. 
\end{proof}
\noindent In Section \ref{3}, we will see examples where the Jacobian is annihilated by a larger ideal in $\ZZ[G]$ than the one generated by the reduced Stickelberger element.

\section{Single Voltage Cyclic Covers of $K_n$ and $K_{n,n}$}
\label{3}

The results in this section were obtained by first computing extensive tables via Sage and Mathematica, from which conjectures were formulated. The conjectures were then proven via lengthy matrix manipulations. See \cite{phdthesis} for additional details. We begin this section by defining the \textit{single voltage assignment}. 

\begin{definition}\label{single voltage}
Let $X$ be a connected graph with vertices labeled as $v_1,\dots,v_n$ such at $v_1$ is adjacent to $v_2$, and let $G$ be the cyclic group of order $d$ generated by $\tau$. Define a \textit{single voltage assignment} to be the voltage assignment $\alpha: E^+(X)\to G$ with $\alpha(e_{1,2})=\tau, \alpha(e_{2,1})=\tau^{-1}$, and whenever there is an edge between $v_i$ and $v_j,$ $\alpha(e_{i,j})=1$ and 0 otherwise.
\end{definition}
%\item Take the complete graph $K_n$, $n>2$ with vertices $v_1,\cdots, v_n$ and put the natural lexicographic resulting orientation on it, namely, $v_i\to v_j$ when $i<j$ (as described above). 
%\item Take a cyclic group of order $d$ generated by $\tau$.
%\item Define \textit{the single voltage assignment} to be the one whose voltage adjacency matrix has $\tau$ in entry $1,2$, $\tau^{-1}$ in entry $2,1,$ and all other off-diagonal entries the identity.
%\end{enumerate}

\noindent Observe that for a single voltage assignment, the resulting voltage adjacency matrix, $A_\alpha$,
which has entries from $\ZZ[G]$ where $G = \langle\tau\rangle$,
``looks almost the same'' as the ordinary adjacency matrix, $A_X$, of $X$. 
Namely, to obtain $A_\alpha$ from $A_X$, simply put $\tau$ in position 1,2 and $\tau^{-1}$ in position 2,1 (instead of the 1's in $A_X$);
and view all other 1's in $A_X$ as the identity of $G$ (in multiplicative notation) to obtain $A_\alpha$.
(Zeros in $A_X$, including along the diagonal, correspondingly remain zero in $A_\alpha.$)\\
\\
\noindent \textit{Remark:} a single voltage assignment depends on the labeling of vertices; 
and different choices of labeling may result in non-isomorphic derived graphs (for example, some connected and some not).
See \cite{phdthesis} for examples.
For a graph whose automorphism group is edge transitive, such as $K_n$, the derived graphs are always isomorphic.

\subsection{Cyclic Covers of $K_n$}

We now let $X=K_n$ in Definition \ref{single voltage}. By symmetry, the Jacobians will \textit{not} depend on which (directed) edge has the non-identity voltage assignment. Putting the single non-identity voltage $\tau$ on edge $v_1\to v_2$ is for convenience. The next proposition follows from \cite{phdthesis}.
\begin{proposition} If $(K_n,Z_d,\alpha)$ is a voltage graph, where $\alpha:E(X)^+\to Z_d$ is the single voltage assignment, then the derived graph $Y$ is connected.
\end{proposition}

\noindent We first determine what the reduced Stickelberger element is of the derived graph $Y$ that corresponds to this voltage graph. We then obtain a $\ZZ[G]$-module presentation of $\Pic(Y).$ Moreover, we show that $\Pic(Y)$ is annihilated by a larger ideal in $\ZZ[G]$ than the one generated by the reduced Stickelberger element. \\

\begin{Theorem}\label{Kn}
Let $(K_n,Z_d,\alpha)$ be as in Definition \ref{single voltage}. 
%Then the voltage Laplacian matrix is the $n \times n$ matrix with $(n{-}1)$'s down the diagonal,
%and all other entries $-1$ except for $-\tau$ in position $1,2$ and $-\tau^{-1}$ in position $2,1$.
Then the reduced Stickelberger element is
$$
\Theta_{Y/X}=-(n-2)  n^{n-3}  (\tau - 1)^2  \tau^{-1}.
$$

\end{Theorem}
\begin{proof}
For $X=K_n$ with single voltage assignment by $Z_d,$ we have the voltage Laplacian matrix with $n-1$ down the diagonal, $-\tau$ in entry $(1,2)$, $-\tau^{-1}$ in entry $(2,1)$ and $-1$ elsewhere:
$$
L_\alpha =  
\begin{pmatrix}
n-1 & -\tau & -1 & -1&\cdots & -1\\
-\tau^{-1}& n-1& -1 &-1&\cdots &-1\\
-1&-1&n-1&-1&\cdots &-1\\
\vdots & \vdots & \vdots &\vdots &\cdots &\vdots\\
-1&-1&-1&-1&\cdots &n-1
\end{pmatrix}
$$
We will now use row and column operations in $\QQ[\tau]$ to put the matrix in essentially upper triangular form. We denote row $i$ and column $j$ of the voltage Laplacian by $C_i$ and $R_j$, respectively. \\
\\
First replace $C_1$ by $C_1+C_3+\cdots +C_n$. %to get
%$$
%\begin{pmatrix}
%1 & -\tau & -1 & -1&\cdots & -1\\
%-\tau^{-1}-n+2& n-1& -1 &-1&\cdots &-1\\
%1&-1&n-1&-1&\cdots &-1\\
%\vdots & \vdots & \vdots &\vdots &\cdots &\vdots\\
%1&-1&-1&-1&\cdots &n-1
%\end{pmatrix}
%$$
Then replace $C_2,\cdots, C_n$ with $C_1+C_2, C_1+C_3,\cdots ,C_1+C_n,$ respectively, to get
$$
\arraycolsep=1.5pt
\left(
\begin{array}{cccccc}
1 & -\tau +1 & 0&0&\cdots & 0\\
-\tau^{-1}-n+2& -\tau^{-1}+1& -\tau^{-1}-n+1 &-\tau^{-1}-n+1&\cdots &-\tau^{-1}-n+1\\
1&0&n&0&\cdots &0\\
\vdots & \vdots & \vdots &\vdots &\cdots &\vdots\\
1&0&0&0&\cdots &n
\end{array}
\right) \qquad (\dagger)
$$
Now replace $C_1$ with $C_1-\frac{1}{n}(C_3+C_4+\cdots+C_n)$ and simplify to obtain
$$
\begin{pmatrix}
1 & -\tau +1 & 0&0&\cdots & 0\\
\frac{-2\tau^{-1}-n+2}{n}& -\tau^{-1}+1& -\tau^{-1}-n+1 &-\tau^{-1}-n+1&\cdots &-\tau^{-1}-n+1\\
0&0&n&0&\cdots &0\\
\vdots & \vdots & \vdots &\vdots &\cdots &\vdots\\
0&0&0&0&\cdots &n
\end{pmatrix}
$$
Let $L'$ be the preceding matrix.
Now use cofactor expansion along the first row of $L'$ to get 
$$
\Theta_{Y/X} = \det L_\alpha = \det L' = 1 \cdot \det L'_{1,1} - (1-\tau) \det L'_{1,2}
$$
where $L'_{i,j}$ is the $i,j$ minor of $L'$.
To compute each of the minor determinants, use cofactor expansion along its first column 
(which has only one nonzero entry), to get:
$$
\det L' = 1 (1 - \tau^{-1}) n^{n - 2} - (1 - \tau) \frac{(-2 \tau^{-1} - n + 2)}{n}  n^{n - 2}.
$$
By elementary algebra this simplifies to the stated conclusion.
\end{proof}
\begin{Corollary}
Under the hypothesis of the Theorem \ref{Kn}, $(n-2)  n^{n-3}  (\tau - 1)^2$ annihilates the Picard group of $Y$, $\Pic(Y)$, hence annihilates its subgroup $\mathcal{J}(Y).$
\end{Corollary}
\begin{proof}
By Corollary \ref{annih} the reduced Stickelberger element $\Theta_{Y/X}$ annihilates the quotient as a $\ZZ[G]$-module. Since $-\tau^{-1}$ is a unit in the ring, the product of the remaining terms effects the annihilation.
\end{proof}

\noindent We will now do further operations on $(\dagger)$ from above to get a $\ZZ[G]$-module presentation of $\Pic(Y).$ We show that $\Pic(Y)$ is annihilated by a larger ideal in $\ZZ[G]$ than the one generated by the reduced Stickelberger element.\\
\\
All of our row and column operations will now be in $\ZZ[\tau]$ (note that to compute the reduced Stickelberger element above, all of the row and column operations were in $\ZZ[\tau]$ until the final step). We begin with the matrix $(\dagger)$.
%$$
%\begin{pmatrix}
%1 & -\tau +1 & 0&0&\cdots & 0\\
%-\tau^{-1}-n+2& -\tau^{-1}+1& -\tau^{-1}-n+1 &-\tau^{-1}-n+1&\cdots &-\tau^{-1}-n+1\\
%1&0&n&0&\cdots &0\\
%\vdots & \vdots & \vdots &\vdots &\cdots &\vdots\\
%1&0&0&0&\cdots &n
%\end{pmatrix}
%$$
Replace $R_2$ with $(R_3+\cdots + R_n)+R_2$.
%$$
%\begin{pmatrix}
%1 & -\tau +1 & 0&0&\cdots & 0\\
%-\tau^{-1}& -\tau^{-1}+1& -\tau^{-1}+1 &-\tau^{-1}+1&\cdots &-\tau^{-1}+1\\
%1&0&n&0&\cdots &0\\
%\vdots & \vdots & \vdots &\vdots &\cdots &\vdots\\
%1&0&0&0&\cdots &n
%\end{pmatrix}
%$$
Then multiply $R_2$ by $-\tau$ and simplify. 
%$$
%\begin{pmatrix}
%1 & -\tau +1 & 0&0&\cdots & 0\\
%1& -\tau+1& -\tau+1 &-\tau+1&\cdots &-\tau+1\\
%1&0&n&0&\cdots &0\\
%\vdots & \vdots & \vdots &\vdots &\cdots &\vdots\\
%1&0&0&0&\cdots &n
%\end{pmatrix}
%$$
Finally replace $R_2$ with $R_2-R_1$ to get
$$
\begin{pmatrix}
1 & -\tau +1 & 0&0&\cdots & 0\\
0& 0& -\tau+1 &-\tau+1&\cdots &-\tau+1\\
1&0&n&0&\cdots &0\\
\vdots & \vdots & \vdots &\vdots &\cdots &\vdots\\
1&0&0&0&\cdots &n
\end{pmatrix}
$$
Denote this matrix by $M$. Then we have
$$
\ZZ^n / L_\alpha(\ZZ^n) = \ZZ^n / M(\ZZ^n).
$$
For $G=Z_d,$ let $e_1, e_2, ...., e_n$ be the standard $\ZZ[G]$-basis vectors. So $e_i$ has a 1 in row $i$ and zeros elsewhere.
Then the relations that define the cokernel of $M$ as a $\ZZ[G]$-module are
 \renum
\item $ e_1 + e_3 + e_4 + .... + e_n = 0$,
\item $(1-\tau) e_1 = 0$,  and
\item $(1-\tau) e_2  +  n e_i = 0,$ \quad   for $3 \le i \le n.$
\end{enumerate}
Now  (i) can be rewritten as
$$
e_1 = - (e_3 + e_4 + ... + e_n)
$$
and so $e_1$ is in the module spanned by $e_3, ...., e_n,$
and therefore, the cokernel of M is generated by just $e_2, e_3, .... , e_n$. So the cokernel of $M$ is the free $\ZZ[G]$-module
generated by $e_2, e_3, ...., e_n$ subject to the relations (via reduction of the above relations):
 \renum
\item $(1-\tau) (e_3 + e_4 + .... + e_n) = 0,$ and
\item $(1-\tau) e_2  +  n e_i  = 0$, \quad  for $3 \le i \le n$
\end{enumerate}
Indeed these give a $\ZZ[G]$-module presentation of $\Pic(Y).$\\
\\
Adding (ii) together for all $i$, we get
$$
(n-2) (1-\tau) e_2 = -n ( e_3 + ... + e_n ).
$$
Multiplying both sides by $(1-\tau)$ and using relation (i), we get
$$
(n-2) (1-\tau)^2 e_2 = n (1 - \tau) (e_3 + ... + e_n) = 0,
$$
hence $(n-2)(1-\tau)^2$ annihilates $e_2.$ For $i \ge 3$ we have 
$$
n e_i = - (1-\tau) e_2.
$$
Multiplying both sides by $(n-2)(1-\tau)$, we get
$$
n (n-2)(1-\tau)e_i = - (n-2)(1-\tau)^2 e_2=0.
$$
So $n(n-2)(1-\tau)$ annihilates $e_i$ for $i=3,\cdots, n.$ This yields the following:
\begin{Theorem}
Under the hypothesis of the Theorem \ref{Kn}, $n(n-2) (\tau - 1)^2$ annihilates $\Pic(Y)$. More specifically, $(n-2)(1-\tau)^2$ annihilates $e_2$ and $n(n-2)(1-\tau)$ annihilates $e_i$ for $i=3,\cdots,d$ where $e_2,\cdots, e_n$ generate the cokernel of $L_\alpha.$
\end{Theorem}

\noindent Using these relations for $\Pic(Y)$, we now compute the Jacobian of $Y$. \\
\\
\noindent We will now assume $n\geq 4$. We will rewrite our matrix relations by changing each $e_i$ to $e_{i-1}$. This yields the following $(n-1)\times (n-1)$ relations matrix
$$
\begin{pmatrix}
0&1-\tau &1-\tau & \cdots & 1-\tau\\
1-\tau & n &0 &\cdots & 0\\
1-\tau & 0 &n & \cdots & 0\\
\vdots &\vdots & \vdots &\ddots &\vdots\\
1-\tau & 0 &0&\cdots &n
\end{pmatrix}
$$
\noindent We will put this matrix in block-diagonal form as follows. First replace $R_3,\cdots, R_{n-1}$ by $R_3-R_2,\cdots, R_{n-1}-R_2,$ respectively.
%$$
%\begin{pmatrix}
%0&1-\tau &1-\tau & \cdots & 1-\tau\\
%1-\tau & n &0 &\cdots & 0\\
%0 & -n &n & \cdots & 0\\
%\vdots &\vdots & \vdots &\ddots &\vdots\\
%0 & -n &0&\cdots &n
%\end{pmatrix}
%$$
Then replace $C_2$ with $C_2+C_3+\cdots+C_{n-1}$,
%$$
%\begin{pmatrix}
%0&(n-2)(1-\tau) &1-\tau & \cdots & 1-\tau\\
%1-\tau & n &0 &\cdots & 0\\
%0 & 0 &n & \cdots & 0\\
%\vdots &\vdots & \vdots &\ddots &\vdots\\
%0 & 0 &0&\cdots &n
%\end{pmatrix}
%$$
and then replace $C_4,\cdots, C_{n-1}$ with $C_4-C_3, \cdots, C_{n-1}-C_3,$ respectively.
%$$
%\begin{pmatrix}
%0&(n-2)(1-\tau) &1-\tau &0& \cdots & 0\\
%1-\tau & n &0 &0&\cdots & 0\\
%0 & 0 &n &-n& \cdots & -n\\
%\vdots &\vdots & \vdots &\vdots&\cdots &\vdots\\
%0 & 0 &0&0&\cdots &n
%\end{pmatrix}
%$$
Lastly, replace $R_3$ with $R_3+\cdots+R_{n-1}$ to get
$$
\begin{pmatrix}
0&(n-2)(1-\tau) &1-\tau &0& \cdots & 0\\
1-\tau & n &0 &0&\cdots & 0\\
0 & 0 &n &0& \cdots & 0\\
\vdots &\vdots & \vdots &\vdots&\cdots &\vdots\\
0 & 0 &0&0&\cdots &n
\end{pmatrix}  \qquad \qquad (*)
$$
From this, we get a $3\times 3$ block matrix in the upper left and a $(n-4)\times (n-4)$ scalar matrix with $n$ on the diagonal on the lower right. Since the lower block is diagonal over $\ZZ[G],$ it remains to further reduce the upper submatrix. Do the following operations on the $3\times 3$ matrix. First replace $C_2$ with $-(n-2)C_3+C_2$. Then
%$$
%\begin{pmatrix}
%0&0 &1-\tau \\
%1-\tau & n &0 \\
%0 & n(2-n) &n
%\end{pmatrix}
%$$
interchange the columns to get 
$$
\begin{pmatrix}
1-\tau&0 &0 \\
0 & 1-\tau &n \\
n & 0 &n(2-n)
\end{pmatrix}
$$
This appears to be the simplest $\ZZ[G]$-module relations matrix for $\text{Pic}(Y).$ \\
\\
Note that doing elementary row and column operations on the voltage Laplacian is the same as
changing generators in the domain and range of the $\ZZ[G]$-module $\text{coker}(L_\alpha)$.
In particular, these also just change generators of its $\ZZ$-module structure.
So in order to compute the Smith Normal Form of the Laplacian of $Y$, we may begin with the already reduced 
voltage Laplacian above. Thus, we will now tensor this matrix with $\rho$ to get $d \times d$ block matrices in each entry, i.e.
$$
\begin{pmatrix}
\rho(1-\tau) &0_d&0_d\\
0_d&\rho(1-\tau)&nI_d\\
nI_d&0_d&n(2-n)I_d
\end{pmatrix}
$$
where $I_d$ denotes the $d\times d$ identity matrix and $0_d$ denotes the $d\times d$ zero matrix. The resulting integer matrix is equivalent to the ``upper left $3\times 3$'' portion of the $\ZZ$-module Laplacian matrix. This yields the following $3d\times 3d$ matrix. 
\setcounter{MaxMatrixCols}{20}
$$
\arraycolsep=3pt
\left(
\begin{array}{ccccccccccccccccccccc}
1&0&\cdots &0&-1&0&0&\cdots &0&0&0&0&\cdots &0\\
-1&1&\cdots &0&0&0&0&\cdots &0&0&0&0&\cdots &0\\
\vdots&\vdots&\cdots &\vdots&\vdots &\vdots&\vdots&\ddots&\vdots&\vdots&\vdots&\vdots&\ddots&\vdots\\
0&0&\cdots &-1&1&0&0&\cdots&0&0&0&0&\cdots&0\\
0&0&\cdots&0&0&1&0&\cdots&0&-1&n&0&\cdots&0\\
0&0&\cdots&0&0&-1&1&\cdots&0&0&0&n&\cdots&0\\
\vdots&\vdots&\ddots &\vdots&\vdots &\vdots&\vdots&\cdots&\vdots&\vdots&\vdots&\vdots&\ddots&\vdots\\
0&0&\cdots &0&0&0&0&\cdots&-1&1&0&0&\cdots&n\\
n&0&\cdots&0&0&0&0&\cdots&0&0&n(2-n)&0&\cdots&0\\
0&n&\cdots &0&0&0&0&\cdots&0&0&0&n(2-n)&\cdots &0\\
\vdots&\vdots&\cdots &\vdots&\vdots &\vdots&\vdots&\ddots&\vdots&\vdots&\vdots&\cdots&\ddots&\vdots\\
0&0&\cdots &0&n&0&0&\cdots &0&0&0&0&\cdots&n(2-n)
\end{array}
\right) 
$$ 
We first do row and column operations that put the $1-\tau$ block in Smith Normal Form, where $\rho(1-\tau)$ is the matrix 
%$$
%\begin{pmatrix}
%1&0&0&\cdots &0&-1\\
%-1&1&0&\cdots &0&0\\
%0&-1&1&\cdots &0&0\\
%\vdots&\vdots&\vdots&\cdots &\vdots&\vdots \\
%0&0&0&\cdots &-1&1
%\end{pmatrix} \qquad  \qquad (**)
%$$
in the upper left corner and middle. \\
\\
For explicit matrices corresponding to the following steps, see \cite{phdthesis}.\\
\\
\textit{Step 1}: Replace $R_1$ and $R_{d+1}$ with $R_1+R_2+\cdots+R_d$ and $R_{d+1}+R_{d+2}+\cdots+R_{2d}$, respectively to get the first row to be all zeros. Then replace $C_1$ and $C_{d+1}$ with $C_1+C_2+\cdots +C_{d}$ and $C_{d+1}+C_{d+2}+\cdots+C_{2d}$, respectively to get the $d+1$-st column to be all zeros. \\
\\
\textit{Step 2}: To zero out all of the $-1'$s on the subdiagonal, replace $R_3,\cdots, R_{d}$ with $R_2+R_3,\cdots, R_{d-1}+R_d$, respectively. Similarly, replace $R_{d+3}, \cdots , R_{2d}$ with\\ $R_{d+2}+R_{d+3},\cdots,R_{2d-1}+R_{2d},$ respectively. This yields
%$$
%\arraycolsep=3pt
%\left(
%\begin{array}{ccccccccccccccccccccc}
%0&0&\cdots &0&0&0&0&\cdots &0&0&0&0&0&\cdots &0\\
%0&1&\cdots &0&0&0&0&\cdots &0&0&0&0&0&\cdots &0\\
%\vdots&\vdots&\cdots &\vdots&\vdots &\vdots&\vdots&\ddots&\vdots&\vdots&\vdots&\vdots&\vdots&\ddots&\vdots\\
%0&0&\cdots &0&1&0&0&\cdots&0&0&0&0&0&\cdots&0\\
%0&0&\cdots&0&0&0&0&\cdots&0&0&n&n&n&\cdots&n\\
%0&0&\cdots&0&0&0&1&\cdots&0&0&0&n&0&\cdots&0\\
%\vdots&\vdots&\ddots &\vdots&\vdots &\vdots&\vdots&\cdots&\vdots&\vdots&\vdots&\vdots&\vdots&\cdots&\vdots\\
%0&0&\cdots &0&0&0&0&\cdots&0&1&0&n&n&\cdots&n\\
%n&0&\cdots&0&0&0&0&\cdots&0&0&n(2-n)&0&0&\cdots&0\\
%n&n&\cdots &0&0&0&0&\cdots&0&0&0&n(2-n)&0&\cdots &0\\
%\vdots&\vdots&\cdots &\vdots&\vdots &\vdots&\vdots&\ddots&\vdots&\vdots&\vdots&\vdots&\vdots &\cdots&\vdots\\
%n&0&\cdots &0&n&0&0&\cdots &0&0&0&0&0&\cdots&n(2-n)
%\end{array}
%\right)  \qquad (2)
%$$
a $(d-1)\times (d-1)$ identity block in the upper left corner and in the middle, and a $(d-1)\times (d-1)$ lower triangular block in the middle right with all entries equal to $n$. \\
\\
\textit{Step 3}: Interchange $C_1$ with $C_{d+1}$.\\
\\
%$$
%\arraycolsep=3pt
%\left(
%\begin{array}{ccccccccccccccccccccc}
%0&0&\cdots &0&0&0&0&\cdots &0&0&0&0&0&\cdots &0\\
%0&1&\cdots &0&0&0&0&\cdots &0&0&0&0&0&\cdots &0\\
%\vdots&\vdots&\cdots &\vdots&\vdots &\vdots&\vdots&\ddots&\vdots&\vdots&\vdots&\vdots&\vdots&\ddots&\vdots\\
%0&0&\cdots &0&1&0&0&\cdots&0&0&0&0&0&\cdots&0\\
%0&0&\cdots&0&0&0&0&\cdots&0&0&n&n&n&\cdots&n\\
%0&0&\cdots&0&0&0&1&\cdots&0&0&0&n&0&\cdots&0\\
%\vdots&\vdots&\ddots &\vdots&\vdots &\vdots&\vdots&\cdots&\vdots&\vdots&\vdots&\vdots&\vdots&\cdots&\vdots\\
%0&0&\cdots &0&0&0&0&\cdots&0&1&0&n&n&\cdots&n\\
%0&0&\cdots&0&0&n&0&\cdots&0&0&n(2-n)&0&0&\cdots&0\\
%0&n&\cdots &0&0&n&0&\cdots&0&0&0&n(2-n)&0&\cdots &0\\
%\vdots&\vdots&\cdots &\vdots&\vdots &\vdots&\vdots&\ddots&\vdots&\vdots&\vdots&\vdots&\vdots &\cdots&\vdots\\
%0&0&\cdots &0&n&n&0&\cdots &0&0&0&0&0&\cdots&n(2-n)
%\end{array}
%\right)  \qquad (3)
%$$
We now remove the first row and column of this matrix, leaving us with a $(3d-1)\times (3d-1)$ matrix. By doing this, we are simply eliminating the invariant factor that is equal to zero.\\
\\
\textit{Step 4}: Zero out the $nI_{d-1}$ matrix in the lower left corner by doing the following: replace $R_{2d+i}$ with $-nR_i+R_{2d+i}$ for all $i=2,\cdots, d$.
%$$
%\arraycolsep=3pt
%\left(
%\begin{array}{cccccccccccccccccccc}
%1&\cdots &0&0&0&0&\cdots &0&0&0&0&0&\cdots &0\\
%\vdots&\cdots &\vdots&\vdots &\vdots&\vdots&\ddots&\vdots&\vdots&\vdots&\vdots&\vdots&\ddots&\vdots\\
%0&\cdots &0&1&0&0&\cdots&0&0&0&0&0&\cdots&0\\
%0&\cdots&0&0&0&0&\cdots&0&0&n&n&n&\cdots&n\\
%0&\cdots&0&0&0&1&\cdots&0&0&0&n&0&\cdots&0\\
%\vdots&\ddots &\vdots&\vdots &\vdots&\vdots&\cdots&\vdots&\vdots&\vdots&\vdots&\vdots&\cdots&\vdots\\
%0&\cdots &0&0&0&0&\cdots&0&1&0&n&n&\cdots&n\\
%0&\cdots&0&0&n&0&\cdots&0&0&n(2-n)&0&0&\cdots&0\\
%0&\cdots &0&0&n&0&\cdots&0&0&0&n(2-n)&0&\cdots &0\\
%\vdots&\ddots &\vdots&\vdots &\vdots&\vdots&\ddots&\vdots&\vdots&\vdots&\vdots&\vdots &\cdots&\vdots\\
%0&\cdots &0&0&n&0&\cdots &0&0&0&0&0&\cdots&n(2-n)
%\end{array}
%\right)  
%$$
This gives a $(d-1)\times (d-1)$ and $2d\times 2d$ block.\\
\\
\textit{Step 5}: We compute the Smith Normal Form of the $2d\times 2d$ block. Since the $(d-1)\times (d-1)$ block is an identity matrix, the invariant factors of the 
full matrix are the same as those of the lower right hand block. \\
%$$
%\begin{pmatrix}
%0&0&0&\cdots&0&n&n&n&\cdots&n\\
%0&1&0&\cdots&0&0&n&0&\cdots&0\\
%0&0&1&\cdots&0&0&n&n&\cdots&0\\
%\vdots &\vdots &\vdots &\ddots&\vdots&\vdots&\vdots&\vdots&\ddots&\vdots\\
%0&0&0&\cdots&1&0&n&n&\cdots&n\\
%n&0&0&\cdots&0&n(2-n)&0&0&\cdots&0\\
%n&0&0&\cdots&0&0&n(2-n)&0&\cdots &0\\
%\vdots &\vdots&\ddots&\vdots&\vdots&\vdots&\vdots&\vdots &\cdots&\vdots\\
%n&0&0&\cdots &0&0&0&0&\cdots&n(2-n)
%\end{pmatrix} \qquad (5)
%$$
\\
\textit{Step 6}:  To zero out the $(d-1)\times (d-1)$ lower triangular submatrix in rows $2$ to $d$ (with all entries equal to $n$), replace $C_{d+j}$ with $-nC_{i}+C_{d+j}$ for all $i=2,\cdots, d$ and for all $j=2,\cdots, i.$ \\
% $$
%\begin{pmatrix}
%0&0&0&\cdots&0&n&n&n&\cdots&n\\
%0&1&0&\cdots&0&0&0&0&\cdots&0\\
%0&0&1&\cdots&0&0&0&0&\cdots&0\\
%\vdots &\vdots &\vdots &\ddots&\vdots&\vdots&\vdots&\vdots&\ddots&\vdots\\
%0&0&0&\cdots&1&0&0&0&\cdots&0\\
%n&0&0&\cdots&0&n(2-n)&0&0&\cdots&0\\
%n&0&0&\cdots&0&0&n(2-n)&0&\cdots &0\\
%\vdots &\vdots&\vdots&\ddots&\vdots&\vdots&\vdots&\vdots &\cdots&\vdots\\
%n&0&0&\cdots &0&0&0&0&\cdots&n(2-n)
%\end{pmatrix} \qquad (6)
%$$
\\
\textit{Step 7}: Move $R_1$ to be the last row and move $C_1$ to be the last column.
% $$
%\begin{pmatrix}
%1&0&\cdots&0&0&0&0&\cdots&0&0\\
%0&1&\cdots&0&0&0&0&\cdots&0&0\\
%\vdots &\vdots &\ddots &\vdots&\vdots&\vdots&\vdots&\ddots&\vdots&\vdots\\
%0&0&\cdots&1&0&0&0&\cdots&0&0\\
%0&0&\cdots&0&n(2-n)&0&0&\cdots&0&n\\
%0&0&\cdots&0&0&n(2-n)&0&\cdots &0&n\\
%\vdots &\vdots&\ddots&\vdots&\vdots&\vdots&\vdots&\cdots &\vdots&\vdots\\
%0&0&\cdots &0&0&0&0&\cdots&n(2-n)&n\\
%0&0&\cdots&0&n&n&n&\cdots&n&0\\
%\end{pmatrix} \qquad (7)
%$$
From this we have a $(d-1)\times (d-1)$ identity block in the upper left and a $(d+1)\times (d+1)$ block in the lower right (all other entries zero). \\
\\
\textit{Step 8}: Remove the $(d-1)\times (d-1)$ identity block for the same reasoning as above.\\
\\
 %$$
%\begin{pmatrix}
%n(2-n)&0&0&\cdots&0&n\\
%0&n(2-n)&0&\cdots &0&n\\
%\vdots&\vdots&\vdots&\cdots&\vdots&\vdots\\
%0&0&0&\cdots&n(2-n)&n\\
%n&n&n&\cdots&n&0\\
%\end{pmatrix} \qquad (4)
%$$
One way to reduce this to Smith Normal Form is as follows. \\
\\
\textit{Step 9}: Replace $R_1$ with $(n-2)R_{d+1}+R_1$.\\
% $$
%\begin{pmatrix}
%0&n(n-2)&n(n-2)&\cdots&n(n-2)&n\\
%0&n(2-n)&0&\cdots &0&n\\
%\vdots&\vdots&\vdots&\cdots&\vdots&\vdots\\
%0&0&0&\cdots&n(2-n)&n\\
%n&n&n&\cdots&n&0\\
%\end{pmatrix} \qquad (9)
%$$
\\
\textit{Step 10}: Replace $R_2,R_3,\cdots, R_d$ with $R_1-R_2,R_1-R_3,\cdots, R_1-R_d$, respectively, 
% $$
%\begin{pmatrix}
%0&n(n-2)&n(n-2)&\cdots&n(n-2)&n\\
%0&2n(n-2)&n(n-2)&\cdots &n(n-2)&0\\
%\vdots&\vdots&\vdots&\cdots&\vdots&\vdots\\
%0&n(n-2)&0&\cdots&2n(n-2)&0\\
%n&n&n&\cdots&n&0\\
%\end{pmatrix} \qquad (5)
%$$
where the $(d-1)\times (d-1)$ sub-matrix in the middle has all diagonal entries equal to $2n(n-2)$ and all other entries equal to $n(n-2).$\\
\\
\textit{Step 11}: Replace $R_{3}$ with $R_{2}+R_3$. Then replace $R_4$ with $R_{3}+R_{4}$. Continue this process until we have replaced $R_{d}$ with $R_{d-1}+R_{d}.$ This yields the matrix
 %$$
%\begin{pmatrix}
%0&n(n-2)&n(n-2)&\cdots&n(n-2)&n\\
%0&2n(n-2)&n(n-2)&\cdots &n(n-2)&0\\
%0&3n(n-2)&3n(n-2)&\cdots &2n(n-2)&0\\
%\vdots&\vdots&\vdots&\cdots&\vdots&\vdots\\
%0&nd(n-2)&nd(n-2)&\cdots&nd(n-2)&0\\
%n&n&n&\cdots&n&0\\
%\end{pmatrix} \qquad (6)
%$$
where the $(d-1)\times (d-1)$ block in the middle has the following property: $R_{i}$ has entry $i(n-2)$ in the first $(i-1)$-st columns and $(i-1)(n-2)$ in the remaining columns, for $i=2,\cdots, d.$\\
\\
\textit{Step 12}: Replace $C_{i}$ for $i=3,\cdots, d$ with $C_{2}-C_{i}$.
% $$
%\begin{pmatrix}
%0&n(n-2)&0&\cdots&0&n\\
%0&2n(n-2)&n(n-2)&\cdots &n(n-2)&0\\
%0&3n(n-2)&0&\cdots &n(n-2)&0\\
%\vdots&\vdots&\vdots&\ddots&\vdots&\vdots\\
%0&n(d-1)(n-2)&0&\cdots&n(n-2)&0\\
%0&nd(n-2)&0&\cdots&0&0\\
%n&n&0&\cdots&0&0\\
%\end{pmatrix} \qquad (12)
%$$
This gives us a $(d-2)\times (d-2)$ upper triangular sub-matrix with all entries equal to $n(n-2).$\\
\\
\textit{Step 13}: Use $C_{3}$ to zero out the remaining $n(n-2)$ entries in row 2, columns $C_{i}$ for $i=4,\cdots, d$ by replacing $C_{i}$ with $C_{i}-C_3.$\\
 %$$
%\begin{pmatrix}
%0&n(n-2)&0&\cdots&0&n\\
%0&2n(n-2)&n(n-2)&\cdots &0&0\\
%0&3n(n-2)&0&\cdots &n(n-2)&0\\
%\vdots&\vdots&\vdots&\ddots&\vdots&\vdots\\
%0&n(d-1)(n-2)&0&\cdots&n(n-2)&0\\
%0&nd(n-2)&0&\cdots&0&0\\
%n&n&0&\cdots&0&0\\
%\end{pmatrix} \qquad (13)
%$$
\\
\textit{Step 14}: Use $C_{4}$ to zero out the remaining $n(n-2)$ in row 3, columns $C_{i}$ for $i=5,\cdots, d$ by replacing $C_{i}$ with $C_{i}-C_{4}.$
 %$$
%\begin{pmatrix}
%0&n(n-2)&0&\cdots&0&n\\
%0&2n(n-2)&n(n-2)&\cdots &0&0\\
%0&3n(n-2)&0&\cdots &0&0\\
%\vdots&\vdots&\vdots&\ddots&\vdots&\vdots\\
%0&n(d-1)(n-2)&0&\cdots&n(n-2)&0\\
%0&nd(n-2)&0&\cdots&0&0\\
%n&n&0&\cdots&0&0\\
%\end{pmatrix} \qquad (14)
%$$
Continue this process until only $n(n-2)$ is left on the diagonal of this $(d-2)\times (d-2)$ block and all other entries are equal to zero.\\
\\
\textit{Step 15}: To zero out all the entries from $R_{2}$ to $R_{d-1}$ in $C_{2}$, first replace $C_{2}$ with $-2C_{3}+C_{2}.$ Then replace $C_{2}$ with $-3C_{4}+C_{2}$. Continue this process until we have
 $$
\begin{pmatrix}
0&n(n-2)&0&\cdots&0&n\\
0&0&n(n-2)&\cdots &0&0\\
0&0&0&\cdots &0&0\\
\vdots&\vdots&\vdots&\ddots&\vdots&\vdots\\
0&0&0&\cdots&n(n-2)&0\\
0&nd(n-2)&0&\cdots&0&0\\
n&n&0&\cdots&0&0\\
\end{pmatrix} 
$$
\textit{Step 16}: Replace $C_{2}$ with $-C_{1}+C_{2}$ and replace $C_{2}$ with $-(n-2)C_{d+1}+C_{2}.$\\
\\
 %$$
%\begin{pmatrix}
%0&0&0&\cdots&0&n\\
%0&0&n(n-2)&\cdots &0&0\\
%0&0&0&\cdots &0&0\\
%\vdots&\vdots&\vdots&\ddots&\vdots&\vdots\\
%0&0&0&\cdots&n(n-2)&0\\
%0&nd(n-2)&0&\cdots&0&0\\
%n&0&0&\cdots&0&0\\
%\end{pmatrix} \qquad (8)
%$$
Permute the columns to obtain a diagonal matrix with the following entries on the diagonal: $n$ with multiplicity 2, $n(n-2)$ with multiplicity $d-2$ and $nd(n-2)$ with multiplicity 1. \\
\\
Finally, returning to the $(n-1)\times (n-1)$ reduced voltage Laplacian $(*)$, when we tensor the $(n-4)\times (n-4)$ matrix in the lower right with $\rho$, 
%$$
%\begin{pmatrix}
%0&(n-2)(1-\tau) &1-\tau &0& \cdots & 0\\
%1-\tau & n &0 &0&\cdots & 0\\
%0 & 0 &n &0& \cdots & 0\\
%\vdots &\vdots & \vdots &\vdots&\cdots &\vdots\\
%0 & 0 &0&0&\cdots &n
%\end{pmatrix}
%$$
we get the entry $n$ on the diagonal with multiplicity $(n-4)d.$ Putting these two matrices together yields the result for $n\geq 4$ and $d\geq 3.$
\begin{Corollary}
For $Y$ as above, the number of spanning trees of $Y$ is
$n^{(n-3)d+1}\cdot (n-2)^{d-1}\cdot d$
\end{Corollary}

\noindent For $d=2,$ apply Steps 1-10. Then apply Step 16. Putting this matrix together with $(*)$ yields the following: 
$$
J(Y) \iso (\ZZ/n\ZZ)^{2(n-4)+2} \oplus (\ZZ/2n(n-2)\ZZ)
$$
 where the exponents indicate the multiplicities of the (distinct) invariant factors. \\
\\
Now when $n=3$ and $d\geq 1$, $X$ is a triangle and $Y$ is a cycle of length $3d.$ This case follows from the familiar Jacobian of a cycle. So we have that
$$
J(Y) \iso \ZZ/3d\ZZ.
$$ 
Alternatively, this can be done independently by corresponding (albeit easier) steps for the $n>3$ reduction above. See \cite{phdthesis} for details.

\subsection{Partial Results on the Jacobian of Single Voltage Cyclic Covers of $K_{n,n}$}\label{3.5}

Now we let $X=K_{n,n}$ in Definition \ref{single voltage}. We partition the $2n$ vertices by: each $v_i$ for $i$ odd is adjacent to every $v_j$ for $j$ even, and vice versa. Note that the Jacobian does not depend on which edge is assigned the nontrivial voltage $\tau$. However, as in Definition \ref{single voltage}, we label the edge $v_1\rightarrow v_2$ with voltage $\tau.$ By \cite{phdthesis}, we get the following:
\begin{proposition} If $(K_{n,n},Z_d,\alpha)$ is a voltage graph, where $\alpha:E(X)^+\to Z_d$ is the single voltage assignment, then the derived graph $Y$ is connected.
\end{proposition}

\noindent We now compute the Smith Normal Form of the Laplacian matrix for the single voltage cyclic cover of $K_{n,n}$ over $\ZZ_{({p})},$ the integers localized at $({p})$, for $p$ not dividing $n.$ This gives us the primes $p$ (and their powers) that divide $|\mathcal{J}(Y)|,$ for $p\nmid n$.\\
\\
We first begin by proving what the reduced Stickelberger element is for the derived graph corresponding to single voltage cyclic covers of $K_{n,n}.$
\begin{Theorem}\label{Kn,n}
Let $(K_{n,n},Z_d,\alpha)$ be as in Definition \ref{single voltage}. %Then the voltage Laplacian matrix is the $2n \times 2n$ matrix with $n$'s down the diagonal, $-\tau$ in entry $(1,2)$, $-\tau^{-1}$ in entry $(2,1)$ and entry $i,j$ equals $-1$ if $i+j$ is odd, and all other entries are zero. 
Then the reduced Stickelberger element is
$$
\Theta_{Y/X}=-(n - 1)^2n^{2n-4}  (\tau - 1)^2  \tau^{-1}.
$$

\end{Theorem}
\begin{proof}
For $X=K_{n,n}$ with single voltage assignment by $Z_d,$ we have the voltage Laplacian matrix with $n$ down the diagonal, $-\tau$ in entry $(1,2)$, $-\tau^{-1}$ in entry $(2,1)$ and entry $i,j$ equals $-1$ if $i+j$ is odd, and all other entries are zero. This gives us the following $2n\times 2n$ matrix:
$$
L_\alpha =  
\begin{pmatrix}
n & -\tau & 0 & -1&0&\cdots &0& -1\\
-\tau^{-1}& n& -1 &0&-1&\cdots &-1&0\\
0&-1&n&-1&0&\cdots &0&-1\\
-1 & 0 & -1 & n&-1&\cdots &-1& 0\\
\vdots & \vdots & \vdots &\vdots &\vdots &\cdots&\vdots&\vdots\\
0&-1&0&-1&0&\cdots &n&-1\\
-1&0&-1&0&-1&\cdots &-1&n
\end{pmatrix}
$$
We will now use row and column operations in $\QQ[\tau]$ to put the matrix in essentially upper triangular form. We denote row $i$ and column $j$ of the voltage Laplacian by $C_i$ and $R_j$, respectively. \\
\\
First replace $C_1$ by $C_1+C_3+\cdots +C_{2n}$.
%$$
%\begin{pmatrix}
%1 & -\tau & 0 & -1&0&\cdots &0& -1\\
%-\tau^{-1}-(n-1)& n& -1 &0&-1&\cdots &-1&0\\
%1&-1&n&-1&0&\cdots &0&-1\\
%0 & 0 & -1 & n&-1&\cdots &-1& 0\\
%\vdots & \vdots & \vdots &\vdots &\vdots &\cdots&\vdots&\vdots\\
%1&-1&0&-1&0&\cdots &n&-1\\
%0&0&-1&0&-1&\cdots &-1&n
%\end{pmatrix}
%$$
Then replace $C_2,C_4,C_6,\cdots, C_{2n}$ with $C_1+C_2, C_1+C_4, C_1+C_6,\cdots, C_1+C_{2n}$, respectively.
 %$$
%\arraycolsep=3pt
%\left(
%\begin{array}{ccccccccccccccccccccc}
%1 & 1-\tau & 0 & 0&0&\cdots &0& 0\\
%-\tau^{-1}-(n-1)& -\tau^{-1}+1& -1 &-\tau^{-1}-(n-1)&-1&\cdots &-1&-\tau^{-1}-(n-1)\\
%1&0&n&0&0&\cdots &0&0\\
%0 & 0 & -1 & n&-1&\cdots &-1& 0\\
%\vdots & \vdots & \vdots &\vdots &\vdots &\cdots&\vdots&\vdots\\
%1&0&0&0&0&\cdots &n&0\\
%0&0&-1&0&-1&\cdots &-1&n
%\end{array}
%\right) 
%$$
This eliminates all of the $-1$'s in all of the even-indexed columns.  We will now eliminate all of the $-1$'s in all of the odd-indexed columns by doing the following: replace $C_3,C_5,C_7,\cdots, C_{2n-1}$ with $\frac{1}{n}C_{2n}+C_3, \frac{1}{n}C_{2n}+C_5, \frac{1}{n}C_{2n}+C_7,\cdots, \frac{1}{n}C_{2n}+C_{2n-1}$.
%$$
%\arraycolsep=1pt
%\begin{pmatrix}
%1 & 1-\tau & 0 & 0&\cdots &0& 0\\
%-\tau^{-1}-(n-1)& -\tau^{-1}+1& \frac{-\tau^{-1}-2n+1}{n} &-\tau^{-1}-(n-1)&\cdots &\frac{-\tau^{-1}-2n+1}{n}&-\tau^{-1}-(n-1)\\
%1&0&n&0&\cdots &0&0\\
%0 & 0 & -1 & n&\cdots &-1& 0\\
%\vdots & \vdots & \vdots &\vdots &\cdots&\vdots&\vdots\\
%1&0&0&0&\cdots &n&0\\
%0&0&0&0&\cdots &0&n
%\end{pmatrix}
%$$
This eliminates the $-1$'s in row $2n$. Now replace $C_3,C_5,C_7,\cdots, C_{2n-1}$ with $\frac{1}{n}C_{2n-2}+C_3, \frac{1}{n}C_{2n-2}+C_5, \frac{1}{n}C_{2n-2}+C_7,\cdots, \frac{1}{n}C_{2n-2}+C_{2n-1}$, respectively.
%$$
%\arraycolsep=1pt
%\begin{pmatrix}
%1 & 1-\tau & 0 & 0&\cdots &0& 0\\
%-\tau^{-1}-(n-1)& -\tau^{-1}+1& \frac{-2\tau^{-1}-3n+2}{n} &-\tau^{-1}-(n-1)&\cdots &\frac{-2\tau^{-1}-3n+2}{n}&-\tau^{-1}-(n-1)\\
%1&0&n&0&\cdots &0&0\\
%0 & 0 & -1 & n&\cdots &-1& 0\\
%\vdots & \vdots & \vdots &\vdots &\cdots&\vdots&\vdots\\
%1&0&0&0&\cdots &n&0\\
%0&0&0&0&\cdots &0&n
%\end{pmatrix}
%$$
This eliminates the $-1$'s in row $2n-2.$ Now replace $C_3,C_5,C_7,\cdots, C_{2n-1}$ with $\frac{1}{n}C_{2n-4}+C_3, \frac{1}{n}C_{2n-4}+C_5, \frac{1}{n}C_{2n-4}+C_7,\cdots, \frac{1}{n}C_{2n-4}+C_{2n-1}$, respectively.
%$$
%\arraycolsep=1pt
%\begin{pmatrix}
%1 & 1-\tau & 0 & 0&\cdots &0& 0\\
%-\tau^{-1}-(n-1)& -\tau^{-1}+1& \frac{-3\tau^{-1}-4n+3}{n} &-\tau^{-1}-(n-1)&\cdots &\frac{-3\tau^{-1}-4n+3}{n}&-\tau^{-1}-(n-1)\\
%1&0&n&0&\cdots &0&0\\
%0 & 0 & -1 & n&\cdots &-1& 0\\
%\vdots & \vdots & \vdots &\vdots &\cdots&\vdots&\vdots\\
%1&0&0&0&\cdots &n&0\\
%0&0&0&0&\cdots &0&n
%\end{pmatrix}
%$$
This eliminates the $-1$'s in row $2n-4.$ Continue this process to zero out the remaining $-1$'s in  $R_{2n-6},\cdots, R_4$ to get the following
$$
\arraycolsep=2pt
\begin{pmatrix}
1 & 1-\tau & 0 & 0&\cdots & 0\\
-\tau^{-1}-(n-1)& -\tau^{-1}+1& \frac{-(n-1)\tau^{-1}-n^2+n-1}{n} &-\tau^{-1}-(n-1)&\cdots &-\tau^{-1}-(n-1)\\
1&0&n&0&\cdots &0\\
0 & 0 & 0 & n&\cdots & 0\\
\vdots & \vdots & \vdots &\vdots &\cdots&\vdots\\
1&0&0&0&\cdots &0\\
0&0&0&0&\cdots &n
\end{pmatrix}
$$
It only remains to zero out the $1$'s in $C_1$ from $R_3$ to $R_{2n-1}$. To do this, we first replace $C_1$ with $-\frac{1}{n}C_{2n-1}+C_1.$ Then replace $C_1$ with $-\frac{1}{n}C_{2n-3}+C_1.$ Continue this process until we replace $C_1$ with $-\frac{1}{n}C_{3}+C_1.$ After simplifying, this yields the following matrix
$$
\arraycolsep=1pt
\begin{pmatrix}
1 & 1-\tau & 0 & 0&\cdots & 0\\
\frac{(-2n+1)\tau^{-1}-(n - 1)^2}{n^2}& -\tau^{-1}+1& \frac{-(n-1)\tau^{-1}-n^2+n-1}{n} &-\tau^{-1}-(n-1)&\cdots &-\tau^{-1}-(n-1)\\
0&0&n&0&\cdots &0\\
0 & 0 & 0 & n&\cdots & 0\\
\vdots & \vdots & \vdots &\vdots &\cdots&\vdots\\
0&0&0&0&\cdots &0\\
0&0&0&0&\cdots &n
\end{pmatrix}
$$
Let $L'$ be the preceding matrix.
Now use cofactor expansion along the first row of $L'$ to get 
$$
\Theta_{Y/X} = \det L_\alpha = \det L' = 1 \cdot \det L'_{1,1} - (1-\tau) \det L'_{1,2}
$$
where $L'_{i,j}$ is the $i,j$ minor of $L'$.
To compute each of the minor determinants, use cofactor expansion along its first column 
(which has only one nonzero entry), to get:
$$
\det L' = 1 (1 - \tau^{-1}) n^{2n - 2} - (1 - \tau) \frac{(-2n+1)\tau^{-1}-(n - 1)^2}{n^2}  n^{2n - 2}.
$$
By elementary algebra this simplifies to the stated conclusion.
\end{proof}

\noindent We will now compute the Smith Normal Form of the previous matrix over $\ZZ_{({p})},$ the integers localized at $({p})$, for $p$ not dividing $n$, so $n$ becomes a unit.\\
\\
Taking the last matrix from above, first replace $R_2$ by $n^2R_2$. Then replace $R_3, \cdots, R_{2n}$ by $\frac{1}{n}R_3,\cdots, \frac{1}{n}R_{2n}$, respectively, to get 
$$
\arraycolsep=8pt
\left(
\begin{array}{ccccc}
1 & 1-\tau & \cdots & 0\\
(-2n+1)\tau^{-1}-(n - 1)^2& -n^2\tau^{-1}+n^2&\cdots &-n^2\tau^{-1}-n^3+n^2\\
0&0&\cdots &0\\
0 & 0 & \cdots & 0\\
\vdots & \vdots & \cdots&\vdots\\
0&0&\cdots &0\\
0&0&\cdots &1
\end{array}
\right)
$$
Now use the $(2n-2)\times (2n-2)$ block identity matrix to zero out the entries in row 2 from columns $3$ through $2n$. \\
\\
Note that doing elementary row and column operations on the voltage Laplacian is the same as
changing generators in the domain and range of the $\ZZ_{({p})}[G]$-module $\text{coker}(L_\alpha)$.
In particular, these also just change generators of its $\ZZ_{({p})}$-module structure.
So in order to compute the Smith Normal Form of the Laplacian of $Y$, we may begin with the already reduced 
voltage Laplacian above. Thus, we will now tensor this matrix with $\rho$ to get the $d \times d$ block matrices in each entry.
$$
\begin{pmatrix}
I_{d} & \rho(1-\tau)\\
(-2n+1)\cdot \rho(\tau^{-1})-(n - 1)^2I_d& -n^2\cdot\rho(\tau^{-1})+n^2\cdot I_d
\end{pmatrix}
$$
where $I_d$ denotes the $d\times d$ identity matrix, 
\begin{align*}
\rho(1-\tau)=&
\begin{pmatrix}
1&0&\cdots &0&-1\\
-1&1&\cdots& 0&0\\
0&-1&\cdots&0&0\\
\vdots &\vdots &\cdots &\vdots &\vdots\\
0&0&\cdots &1&0\\
0&0&\cdots &-1&1
\end{pmatrix}\\
\\
-n^2\cdot\rho(\tau^{-1})+n^2\cdot I_d=&
\begin{pmatrix}
n^2&-n^2&\cdots &0&0\\
0&n^2&\cdots &0&0\\
0&0&\cdots &0&0\\
\vdots &\vdots &\cdots &\vdots &\vdots\\
0&0&\cdots &0&-n^2\\
-n^2&0&\cdots &0&n^2
\end{pmatrix}\\
\end{align*}
and lastly,
$$
(-2n+1)\cdot \rho(\tau^{-1})-(n - 1)^2I_d=
\begin{pmatrix}
-(n - 1)^2&-2n+1&\cdots&0&0\\
0&-(n - 1)^2&\cdots &0&0\\
\vdots&\vdots &\ddots &\vdots &\vdots\\
0&0&\cdots & -(n - 1)^2&-2n+1\\
-2n+1&0&\cdots &0&-(n - 1)^2
\end{pmatrix}
$$
Putting these together, we get the following matrix
$$
\begin{pmatrix}
1&0&\cdots &0 &1&0&\cdots &0&-1\\
0&1&\cdots &0&-1&1&\cdots &0&0\\
\vdots &\vdots &\ddots &\vdots &\vdots &\vdots &\cdots &\vdots &\vdots\\
0&0&\cdots&1&0&0&\cdots &-1&1\\
-(n - 1)^2&-2n+1&\cdots &0&n^2&-n^2&\cdots &0&0\\
0&-(n - 1)^2&\cdots &0&0&n^2&\cdots &0&0\\
\vdots &\vdots &\cdots &\vdots &\vdots &\vdots &\ddots &\vdots &\vdots\\
0&0&\cdots &-2n+1&0&0&\cdots &n^2&-n^2\\
-2n+1&0&\cdots &-(n - 1)^2&-n^2&0&\cdots &0&n^2
\end{pmatrix}
$$
We will now use row and column operations to put this matrix in diagonal form. To see explicit matrices of the steps that follow, refer to \cite{phdthesis}.\\
\\
\textit{Step 1:} To zero out the last column, replace $C_{2d}$ with $C_{d+1}+C_{d+2}+\cdots C_{2d}$.\\
\\
%$$
%\begin{pmatrix}
%1&0&\cdots &0 &1&0&\cdots &0&0\\
%0&1&\cdots &0&-1&1&\cdots &0&0\\
%\vdots &\vdots &\ddots &\vdots &\vdots &\vdots &\cdots &\vdots &\vdots\\
%0&0&\cdots&1&0&0&\cdots &-1&0\\
%-(n - 1)^2&-2n+1&\cdots &0&n^2&-n^2&\cdots &0&0\\
%0&-(n - 1)^2&\cdots &0&0&n^2&\cdots &0&0\\
%\vdots &\vdots &\cdots &\vdots &\vdots &\vdots &\ddots &\vdots &\vdots\\
%0&0&\cdots &-2n+1&0&0&\cdots &n^2&0\\
%-2n+1&0&\cdots &-(n - 1)^2&-n^2&0&\cdots &0&0
%\end{pmatrix}
%$$
\textit{Step 2:} Replace $R_{2d}$ with $R_{d+1}+R_{d+2}+\cdots+R_{2d}$.\\
\\
%$$
%\begin{pmatrix}
%1&0&\cdots &0 &1&0&\cdots &0&0\\
%0&1&\cdots &0&-1&1&\cdots &0&0\\
%\vdots &\vdots &\ddots &\vdots &\vdots &\vdots &\cdots &\vdots &\vdots\\
%0&0&\cdots&1&0&0&\cdots &-1&0\\
%-(n - 1)^2&-2n+1&\cdots &0&n^2&-n^2&\cdots &0&0\\
%0&-(n - 1)^2&\cdots &0&0&n^2&\cdots &0&0\\
%\vdots &\vdots &\cdots &\vdots &\vdots &\vdots &\ddots &\vdots &\vdots\\
%0&0&\cdots &-2n+1&0&0&\cdots &n^2&0\\
%-n^2&-n^2&\cdots &-n^2&0&0&\cdots &0&0
%\end{pmatrix}
%$$
This gives us a matrix with entry $-n^2$ in $R_{2d}, \ C_1$ through $C_{d}$ and $0$ from $C_{d+1}$ through $C_{2d}$.\\
\\
\textit{Step 3:} Replace $R_1$ with $R_1+R_2+\cdots+R_d$.
%$$
%\begin{pmatrix}
%1&1&\cdots &1 &0&0&\cdots &0&0\\
%0&1&\cdots &0&-1&1&\cdots &0&0\\
%\vdots &\vdots &\ddots &\vdots &\vdots &\vdots &\cdots &\vdots &\vdots\\
%0&0&\cdots&1&0&0&\cdots &-1&0\\
%-(n - 1)^2&-2n+1&\cdots &0&n^2&-n^2&\cdots &0&0\\
%0&-(n - 1)^2&\cdots &0&0&n^2&\cdots &0&0\\
%\vdots &\vdots &\cdots &\vdots &\vdots &\vdots &\ddots &\vdots &\vdots\\
%0&0&\cdots &-2n+1&0&0&\cdots &n^2&0\\
%-n^2&-n^2&\cdots &-n^2&0&0&\cdots &0&0
%\end{pmatrix}
%$$
This gives us a matrix with all entries equal to 1 in $R_1, C_{1}$ through $C_{d}$ and $0$ from $C_{d}$ to $C_{2d}.$\\
\\
\textit{Step 4:} To zero out the $1$'s in $R_1,$ $C_2$ through $C_{d}$ and to zero out the $-n^2$ in $R_{2d}$ as well as columns $C_{2}$ through $C_{d},$ replace $C_2, C_3,\cdots, C_d$ with $C_2-C_1, C_3-C_1,\cdots, C_d-C_1,$ respectively. 
%$$
%\begin{pmatrix}
%1&0&0&\cdots &0 &0&0&\cdots &0&0\\
%0&1&0&\cdots &0&-1&1&\cdots &0&0\\
%\vdots &\vdots &\vdots&\cdots &\vdots &\vdots &\vdots &\cdots &\vdots &\vdots\\
%0&0&0&\cdots&1&0&0&\cdots &-1&0\\
%-(n - 1)^2&n^2-4n+2&(n-1)^2&\cdots &(n-1)^2&n^2&-n^2&\cdots &0&0\\
%0&-(n - 1)^2&-2n+1&\cdots &0&0&n^2&\cdots &0&0\\
%\vdots &\vdots &\vdots &\ddots &\vdots &\vdots &\vdots &\ddots &\vdots &\vdots\\
%0&0&0&\cdots &-2n+1&0&0&\cdots &n^2&0\\
%-n^2&0&0&\cdots &0&0&0&\cdots &0&0
%\end{pmatrix}
%$$
In $R_d,$ we have $-(n - 1)^2$ in the first column, $n^2-4n+2$ in the second column, and $(n-1)^2$ in columns $3$ through $d$.\\
\\
\textit{Step 5:} Use $R_1$ to zero out the remaining entries in $C_1$.\\
\\
%$$
%\begin{pmatrix}
%1&0&0&\cdots &0 &0&0&\cdots &0&0\\
%0&1&0&\cdots &0&-1&1&\cdots &0&0\\
%\vdots &\vdots &\vdots&\cdots &\vdots &\vdots &\vdots &\cdots &\vdots &\vdots\\
%0&0&0&\cdots&1&0&0&\cdots &-1&0\\
%0&n^2-4n+2&(n-1)^2&\cdots &(n-1)^2&n^2&-n^2&\cdots &0&0\\
%0&-(n - 1)^2&-2n+1&\cdots &0&0&n^2&\cdots &0&0\\
%\vdots &\vdots &\vdots &\ddots &\vdots &\vdots &\vdots &\ddots &\vdots &\vdots\\
%0&0&0&\cdots &-2n+1&0&0&\cdots &n^2&0\\
%0&0&0&\cdots &0&0&0&\cdots &0&0
%\end{pmatrix}
%$$
\textit{Step 6:} Remove $C_{2d}$ and $R_{2d}$, as they do not contribute any nonzero invariant factors. This leaves us with a $({2d-1})\times ({2d-1})$ matrix:
$$
\begin{pmatrix}
1&0&0&\cdots &0 &0&0&\cdots &0\\
0&1&0&\cdots &0&-1&1&\cdots &0\\
\vdots &\vdots &\vdots&\cdots &\vdots &\vdots &\vdots &\cdots &\vdots \\
0&0&0&\cdots&1&0&0&\cdots &-1\\
0&n^2-4n+2&(n-1)^2&\cdots &(n-1)^2&n^2&-n^2&\cdots &0\\
0&-(n - 1)^2&-2n+1&\cdots &0&0&n^2&\cdots &0\\
\vdots &\vdots &\vdots &\ddots &\vdots &\vdots &\vdots &\ddots &\vdots \\
0&0&0&\cdots &-2n+1&0&0&\cdots &n^2\\
\end{pmatrix}
$$
\textit{Step 7:} Replace $C_{2d-1}$ with $C_{d+1}+C_{d+2}+\cdots+C_{2d-1}$.
%$$
%\begin{pmatrix}
%1&0&0&\cdots &0 &0&0&\cdots &0\\
%0&1&0&\cdots &0&-1&1&\cdots &0\\
%\vdots &\vdots &\vdots&\cdots &\vdots &\vdots &\vdots &\cdots &\vdots \\
%0&0&0&\cdots&1&0&0&\cdots &-1\\
%0&n^2-4n+2&(n-1)^2&\cdots &(n-1)^2&n^2&-n^2&\cdots &0\\
%0&-(n - 1)^2&-2n+1&\cdots &0&0&n^2&\cdots &0\\
%\vdots &\vdots &\vdots &\ddots &\vdots &\vdots &\vdots &\ddots &\vdots \\
%0&0&0&\cdots &-2n+1&0&0&\cdots &n^2\\
%\end{pmatrix}
%$$
This makes all entries in $C_{2d-1}$ zero except for the entries in $R_{d}$ and $R_{2d-1},$ which are $-1$ and $n^2$, respectively. \\
\\
\textit{Step 8:} Continue this process; replace $C_{2d-2}, \cdots, C_{d+2}$ with $C_{d+1}+C_{d+2}+\cdots +C_{2d-2}, \cdots, C_{d+1}+C_{d+2},$ respectively.
%$$
%\begin{pmatrix}
%1&0&0&\cdots &0 &0&0&\cdots &0\\
%0&1&0&\cdots &0&-1&0&\cdots &0\\
%\vdots &\vdots &\vdots&\cdots &\vdots &\vdots &\vdots &\cdots &\vdots \\
%0&0&0&\cdots&1&0&0&\cdots &-1\\
%0&n^2-4n+2&(n-1)^2&\cdots &(n-1)^2&n^2&0&\cdots &0\\
%0&-(n - 1)^2&-2n+1&\cdots &0&0&n^2&\cdots &0\\
%\vdots &\vdots &\vdots &\ddots &\vdots &\vdots &\vdots &\ddots &\vdots \\
%0&0&0&\cdots &-2n+1&0&0&\cdots &n^2\\
%\end{pmatrix}
%$$
This gives us a $(d-1)\times (d-1)$ diagonal block in the upper right with $-1$'s down the diagonal and a $(d-1)\times (d-1)$ diagonal block in the lower right with $n^2$'s down the diagonal.\\
\\
\textit{Step 9:} To zero out the $(d-1)\times (d-1)$ block in the upper right with $-1$'s down the diagonal, replace $C_{d+1}, C_{d+2},\cdots, C_{2d-1}$ with $C_2+C_{d+1}, C_3+C_{d+2},\cdots, C_d+C_{2d-1}$.
%$$
%\begin{pmatrix}
%1&0&0&\cdots &0 &0&0&\cdots &0\\
%0&1&0&\cdots &0&0&0&\cdots &0\\
%\vdots &\vdots &\vdots&\cdots &\vdots &\vdots &\vdots &\ddots &\vdots \\
%0&0&0&\cdots&1&0&0&\cdots &0\\
%0&n^2-4n+2&(n-1)^2&\cdots &0&2(n-1)^2&(n-1)^2&\cdots &(n-1)^2\\
%0&-(n - 1)^2&-2n+1&\cdots &0&-(n - 1)^2&(n-1)^2&\cdots &0\\
%\vdots &\vdots &\vdots &\ddots &\vdots &\vdots &\vdots &\ddots &\vdots \\
%0&0&0&\cdots &-2n+1&0&0&\cdots &(n-1)^2\\
%\end{pmatrix}
%$$
In the lower right $(d-1)\times (d-1)$ block, we have entry $2(n-1)^2$ in column and row $d+1,$ and then entry $(n-1)^2$ in row $d+1,$ columns $d+2$ through $2d-1.$ Then in columns $d+1$ through $2d-1,$ row $d+2$ through row $2d-1$, we have $-(n - 1)^2$ on the lower subdiagonal, and $(n-1)^2$ on the diagonal.\\
\\
\textit{Step 10:} Use the identity $d\times d$ block in the upper left to zero out the $(d-1)\times (d-1)$ block in the lower left. \\
\\
%$$
%\begin{pmatrix}
%1&0&0&\cdots &0 &0&0&\cdots &0\\
%0&1&0&\cdots &0&0&0&\cdots &0\\
%\vdots &\vdots &\vdots&\cdots &\vdots &\vdots &\vdots &\ddots &\vdots \\
%0&0&0&\cdots&1&0&0&\cdots &0\\
%0&0&0&\cdots &0&2 (n - 1)^2&(n-1)^2&\cdots &(n-1)^2\\
%0&0&0&\cdots &0&-(n - 1)^2&(n-1)^2&\cdots &0\\
%\vdots &\vdots &\vdots &\ddots &\vdots &\vdots &\vdots &\ddots &\vdots \\
%0&0&0&\cdots &0&0&0&\cdots &(n-1)^2\\
%\end{pmatrix}
%$$
\textit{Step 11:} Since the $d\times d$ block in the upper left contributes no nontrivial factors, we will only consider the $(d-1)\times (d-1)$ block in the lower right.\\
%$$
%\begin{pmatrix}
%2 (n - 1)^2&(n - 1)^2&(n - 1)^2&\cdots&(n - 1)^2 &(n - 1)^2\\
%-(n - 1)^2& (n - 1)^2&0&\cdots &0&0\\
%0&-(n - 1)^2&(n - 1)^2&\cdots &0&0\\
%\vdots&\vdots&\vdots&\ddots&\vdots&\vdots\\
%0&0&0&\cdots&(n - 1)^2&0\\
%0&0&0&\cdots&-(n - 1)^2&(n - 1)^2
%\end{pmatrix}
%$$
\\
\textit{Step 12:} Replace $C_{d-1}$ with $C_1+C_2+\cdots+C_{d-1}$.\\
\\
%$$
%\begin{pmatrix}
%2 (n - 1)^2&(n - 1)^2&(n - 1)^2&\cdots&(n - 1)^2&d (n - 1)^2\\
%-(n - 1)^2& (n - 1)^2&0&\cdots &0&0\\
%0&-(n - 1)^2&(n - 1)^2&\cdots &0&0\\
%\vdots&\vdots&\vdots&\ddots&\vdots&\vdots\\
%0&0&0&\cdots&(n - 1)^2&0\\
%0&0&0&\cdots&-(n - 1)^2&0
%\end{pmatrix}
%$$
\textit{Step 13:} Replace $C_{d-2}$ with $C_1+C_2+\cdots+C_{d-2}$.\\
\\
%$$
%\begin{pmatrix}
%2 (n - 1)^2&(n - 1)^2&\cdots&(d - 1) (n - 1)^2&d (n - 1)^2\\
%-(n - 1)^2& n^2+2n-1&\cdots &0&0\\
%0&-(n - 1)^2&\cdots &0&0\\
%\vdots&\vdots&\cdots&\vdots&\vdots\\
%0&0&\cdots&0&0\\
%0&0&\cdots&-(n - 1)^2&0
%\end{pmatrix}
%$$
\textit{Step 14:} Continue this process to zero out the remaining $(n-1)^2$ terms on the diagonal. \\
\\
%$$
%\begin{pmatrix}
%2 (n - 1)^2&3(n - 1)^2&\cdots&(d - 1) (n - 1)^2&d (n - 1)^2\\
%-(n - 1)^2& 0&\cdots &0&0\\
%0&-(n - 1)^2&\cdots &0&0\\
%\vdots&\vdots&\cdots&\vdots&\vdots\\
%0&0&\cdots&0&0\\
%0&0&\cdots&-(n - 1)^2&0
%\end{pmatrix}
%$$
\textit{Step 15:} Use the $(d-2)\times (d-2)$ block diagonal matrix with entries $-(n - 1)^2$ to zero out the entries in row 1, columns $1$ through $d-2$.\\
\\
%$$
%\begin{pmatrix}
%0&0&\cdots&0&d (n - 1)^2\\
%-(n - 1)^2& 0&\cdots &0&0\\
%0&-(n - 1)^2&\cdots &0&0\\
%\vdots&\vdots&\cdots&\vdots&\vdots\\
%0&0&\cdots&0&0\\
%0&0&\cdots&-(n -1)^2&0
%\end{pmatrix}
%$$
\textit{Step 16:} Multiply each of $C_1,\cdots, C_{d-1}$ by the unit $-1$ and also move $C_{d-1}$ to be the first column
to get the diagonal matrix:\\
$$
\begin{pmatrix}
d(n-1)^2&0&0&\cdots&0\\
0&(n - 1)^2& 0&\cdots &0\\
0&0&(n - 1)^2&\cdots &0\\
\vdots&\vdots&\vdots&\cdots&\vdots\\
0&0&0&\cdots&0\\
0&0&0&\cdots&(n -1)^2
\end{pmatrix}
$$
%Lastly, replace $C_1,\cdots, C_{d-2}$ by $(-1)C_1, \cdots, (-1)C_{d-2}$ to get
%$$
%\begin{pmatrix}
%0&0&\cdots&0&d (n - 1)^2\\
%(n - 1)^2& 0&\cdots &0&0\\
%0&(n - 1)^2&\cdots &0&0\\
%\vdots&\vdots&\cdots&\vdots&\vdots\\
%0&0&\cdots&0&0\\
%0&0&\cdots&(n - 1)^2&0
%\end{pmatrix}
%$$
By \cite{Ger76}, this gives us the primes $p$ (and their powers) that divide $|\mathcal{J}(Y)|,$ for $p\nmid n$. We write this as the following theorem
\begin{Theorem}\label{partialKnn}
Let $Y$ be a single voltage cover of the complete bipartite graph $K_{n,n}$ by the cyclic group of order $d$, where $n,d\geq 3$. For any prime $p$ with $p\nmid n$, the Sylow p-subgroup, $\mathcal{J}_p(Y)$, of $
\mathcal{J}(Y)$ has the following (elementary divisor) decomposition:
$$
\mathcal{J}_p(Y)\cong (\ZZ/p^{2a}\ZZ)^{d-2}\oplus (\ZZ/p^{2a+b}\ZZ)^1
$$
where $p^a$ is the largest power of $p$ dividing $n-1$ and $p^b$ is the largest power of $p$ dividing $d$. In particular, $|\mathcal{J}_p(Y)|=p^{2a(d-1)+b},$ and the $p$-rank of $\mathcal{J}(Y)$ is $d-1.$
\end{Theorem}

\section{Zeta Functions of Voltage Graphs}\label{4}

Assume all base graphs $X$ are finite and connected. 
We now introduce some terminology. These can also be found in \cite{HMSV19}.
\begin{definition}
A closed path 
$$
P: w=w_1\xrightarrow{e_{1,2}} w_2\xrightarrow{e_{2,3}} \cdots \xrightarrow{e_{m-1,m}} w_m=w
$$
is called a prime path if it has no backtrack or tail and one may only go around the path once. For the closed path $P$, the equivalence class $[P]$ means the following
$$
[P]=\{e_{1,2}e_{2,3}\cdots e_{m-1,m}, \quad e_{2,3}\cdots e_{m-1,m}e_{1,2}, \qquad \cdots  \qquad ,e_{m-1,m}e_{1,2},\cdots e_{m-2,m-1}\}.
$$
That is, two paths are equivalent if we can get one from the other by changing the starting vertex. A prime in a graph $X$ is an equivalence class $[P]$ of prime paths. The length of $P$ is $\nu({P})=m$, the number of edges in $P$.
\end{definition}

\begin{definition}
The Ihara zeta function of a finite connected graph $X$ is defined to be
$$
\zeta_X(u)=\prod_{[P]}(1-u^{\nu({P})})^{-1}
$$
where the product is over all primes $[P]$ in $X$, \ $u\in \CC$ and $|u|$ is sufficiently small.
\end{definition}
\noindent The following theorem, which can be found in \cite{Ter11}, gives a formula for computing the zeta function.
\begin{Theorem}[Three-term determinant formula]{}
Let $A$ be the adjacency matrix of $X$ and let $D$ be the diagonal matrix of vertex degrees (which is a scalar matrix if $X$ is a regular graph). Let
$$
r_X=r=|E(X)|-|V(X)|+1
$$
where $E(X)$ and $V(X)$ denote the edges and vertices of $X$, respectively. Then we have the Ihara three-term determinant formula
$$
\zeta_X(u)^{-1}=(1-u^2)^{r-1}\det(I-Au+(D-I)u^2).
$$
\end{Theorem}

\subsection{Zeta Functions and the Order of the Jacobian}\label{4.1}

Throughout this section, for $(X,G,\alpha)$ a voltage graph, assume
\renum
\item $X$ is connected,
\item $r_X-1\neq 0,$
\item $G$ is abelian, and
\item $\alpha:E(X)^+\to G$ is such that $Y$ is connected (see \cite{phdthesis}).
\end{enumerate}
By Theorem 2.11 in \cite{HMSV19}, we have that 
$$
\zeta^*_X(1)=(-1)^{r+1}2^r(r-1)\kappa_X, 
$$
where $\zeta_X^*(1)$ denotes the first non-vanishing Taylor coefficient of $\zeta_X(u)^{-1}$ at $u=1$ and $\kappa_X$ is the number of spanning trees in $X$. Thus, rearranging
\begin{equation}\label{equation 4.1}
\zeta^*_X(1)/((-1)^{r+1}2^r(r-1))=\kappa_X=|\mathcal{J}(X)|.
\end{equation}
We now seek a formula that relates $|\mathcal{J}(X)|$ to $|\mathcal{J}(Y)|$. The zeta function of a derived graph $Y$ of a connected base graph $X$ with voltage assignment by $G$, where $G$ is abelian, is by Theorem 3.1 in \cite{HMSV19},
\begin{equation}\label{equation 4.2}
\zeta_Y(u)=\zeta_X(u)\prod_{\chi_t\neq \chi_0}L(u,\chi_t)
\end{equation}
where 
$$
L(u,\chi_t)^{-1}=(1-u^2)^{r-1}\det(I-A_{\chi_t} u+(D-I)u^2)
$$ 
with the product of the $\chi_t$ over all the irreducible characters of $G$ and where $\chi_0$ is the trivial character, and $A_{\chi}$ (called the \textit{Artinized adjacency matrix}) is defined to be
$$
A_{\chi}=\sum_{\sigma\in G} \chi(\sigma)A(\sigma)
$$
where (because $X$ has no loops or multiple edges) the $i,j$-entry of $A(\sigma)$ simplifies to
$$
a_{i,j}(\sigma)=\begin{cases}
1& \text{ if } v_{i,1}\rightarrow v_{j,\sigma} \text{ in } Y\\
0 & \text{ otherwise}
\end{cases}
$$
where $1=\tau_0$ is the identity of $G$. \\
\\
Since $Y$ is the derived graph, $a_{i,j}(\sigma)$ is nonzero if and only if $v_i\xrightarrow{\tau}v_j$ in $X$ and $\sigma=\tau$ (in which case $v_j\xrightarrow{\tau^{-1}}v_i$ and $a_{j,i}(\sigma^{-1})\neq 0$ with $\sigma^{-1}=\tau^{-1}$). In other words, the $i,j$-entry of $A(\sigma)$ is nonzero for a unique $\sigma$ in $G$, namely when $\sigma$ equals the voltage of the edge $v_i\to v_j.$ Thus
$$
\sum_{\sigma\in G} \sigma A(\sigma)=A_\alpha \qquad \text{(the voltage adjacency matrix)}.
$$
Then for every irreducible character $\chi$ of $G$
$$
A_\chi=\chi(A_\alpha)
$$
where $\chi(A_\alpha)$ means evaluating every group element $\sigma$ of the voltage adjacency matrix at $\chi$ (and $\chi(0)=0)$.\\
\\
Now using $A_\alpha$, the voltage adjacency matrix, in place of $A_{\chi_t}$ in $L(u,\chi_t)^{-1},$ we define 
\begin{equation}\label{equation 4.6}
L(u,\alpha)^{-1}=(1-u^2)^{r-1}\det(I-A_\alpha u+(D-I)u^2)
\end{equation}
to be the \textit{equivariant $L$-function}. From this, we now define the \textit{reduced equivariant $L$ function} 
\begin{equation}\label{equation 4.3}
L^*(u,\alpha)^{-1}=\det(I-A_\alpha u +(D-I)u^2).
\end{equation}
Then evaluate each group element entry in $A_\alpha$ at the degree-one character $\chi_t$ 
(the zero entries remain zero).
When this matrix is substituted in place of $A_\alpha$ in (\ref{equation 4.6})---since such evaluation is a ring homomorphism from $\ZZ[G]$ to $\CC$---the resulting function is $L(u,\chi_t)^{-1}$; 
and likewise when we substitute the same matrix for $A_\alpha$ in (\ref{equation 4.3}) we get $L^*(u,\chi_t)^{-1}$.\\
\\
Now observe that the equivariant $L$-function can be factored as 
$$
L(u,\alpha)^{-1}=(1 - u)^{r_X - 1}(1 + u)^{r_X - 1}\det(I - A_\alpha u + (D - I)u^2).
$$
Expanding each of these terms of $(u-1)$, we get
\begin{align*}
(1 - u)^{r_X - 1}&=(-1)^{r_X - 1}(u - 1)^{r_X - 1} \qquad \text{and}\\
\\
(1 + u)^{r_X - 1} &= (2 + (u - 1))^{r_X- 1} \\
&= 2^{r_X - 1} + 2^{ r_X - 2}(r_X - 1)(u - 1) +\cdots + (\text{higher powers of $(u-1)$)}.
\end{align*}
The reduced equivariant $L$-function is of the form
$$
c_0+c_1(u-1)+c_2(u-1)^2+\cdots \text{higher powers of $(u-1)$}
$$
where the coefficients $c_i$ are from the integral group ring $\ZZ[G].$\\
\\
So multiplying the above polynomials together, we get the lowest nonzero coefficient (which corresponds to the $(u-1)^{r_X-1}$ factor, which we know is nonzero) to be 
$$
(-1)^{r_X- 1}2^{r_X - 1} c_0.
$$
For $G\neq 1$, letting $u=1$ in formula (\ref{equation 4.3}) for $L^*(u,\alpha)^{-1},$ we see that $c_0$ is exactly the reduced Stickelberger element, $\Theta_{Y/X}$.
We then evaluate the degree-one irreducible characters $\chi_t$ at the group element $\Theta_{Y/X}$, which we denote by $\chi_t(\Theta_{Y/X}),$ and take the product to get
$$
\prod_{\chi_t\neq \chi_0}L(1,\chi_t)=\prod_{\chi_t\neq \chi_0}(-1)^{r_X- 1}2^{r_X - 1}\chi_t(\Theta_{Y/X}).
$$
\\
Lastly, in (\ref{equation 4.2}) we multiply by the coefficient corresponding to the $L$-function for $X$, i.e., $\zeta_X$, which equals $(-1)^{(r_X+1)}2^{r_X}(r_X-1)|\mathcal{J}(X)|$ to get that the first non-vanishing  Taylor coefficient of $\zeta_Y(y)^{-1}$ at $u=1$ is
\begin{equation}\label{equation 4.4}
\zeta^*_Y(1)=(-1)^{r_X+1}2^{r_X}(r_X-1)|\mathcal{J}(X)|\prod_{\chi_t\neq \chi_0}(-1)^{r_X- 1}2^{r_X - 1}\chi_t(\Theta_{Y/X}).
\end{equation}
Now note that for any voltage cover $Y$ of $X$ of degree $d$, we have $r_Y-1=d(r_X-1).$ Also, the number of irreducible characters of $G$ (abelian) is $d$. By formula $(\ref{equation 4.1})$ applied to $Y$, instead of $X$, we have
$$
|\mathcal{J}(Y)|=\zeta_Y^*(1)/((-1)^{r_Y-1}2^{r_Y}(r_Y-1)).
$$
Substituting the value of $\zeta_Y^*(1)$ from (\ref{equation 4.4}) into this and simplifying gives Theorem \ref{OrdJac}.

\subsubsection{Zeta Functions of Cyclic Graph Coverings}\label{4.1.1}

	Assume throughout this subsection that $(X,G,\alpha)$ is a voltage graph such that
\renum
\item $X$ is connected, 
\item $G=\langle \tau\rangle$ is cyclic of order $d$,
\item there is a closed walk in $X$ with net voltage $\tau$, and
\item $r_X-1\neq 0$.
\end{enumerate}
By \cite{phdthesis}, these first three conditions ensure that all covering graphs by $Z_d=\langle \tau \rangle$ are connected.\\
\\
Fix the primitive $d^{th}$ root of unity $\lambda=e^{2\pi i /d}.$ Then every irreducible character of $Z_d$ maps $\tau$ to a $d$th root of unity, i.e.,
$$
\chi_t:\tau\to \lambda^t, \qquad 0\leq t\leq d-1,
$$
where $\chi_0$ is again the principal character of $G$.
\begin{Corollary}\label{OrdJac1}
Let $X$ satisfy the hypotheses at the outset of this section, with voltage assignment by $Z_d.$ Then  the order of the Jacobian of the derived graph $Y$ is
$$
|\mathcal{J}(Y)|=\frac{1}{d}\cdot |\mathcal{J}(X)|\prod_{t=1}^{d-1}\chi_t(\Theta_{Y/X}),
$$
where $\chi_t:\tau\mapsto \lambda^t$. 
\end{Corollary}
\begin{proof}
This follows immediately from Theorem \ref{OrdJac}.
\end{proof}

\noindent We now go on to show that when $(X,Z_d,\alpha)$ is a voltage graph with $\alpha$ given by the \textit{single} voltage assignment, as in Definition \ref{single voltage}, then the reduced Stickelberger element, $\Theta_{Y/X}$, is of a specific form. We first prove the following lemma.
\begin{Lemma}\label{DetL}
Let $L$ be any symmetric $n \times n$ matrix with entries from $\ZZ$ with $n \ge 3$, and let $x$ be an indeterminate over $\ZZ$.
Assume $L$ has 1 in position 1,2 (hence also in position 2,1) and $\det L = 0$.
Let $L_x$ be the same matrix $L$ except with the 1,2 entry replaced by $x$ and the 2,1 entry replaced by $1/x$,
and let $\Theta(x) = \det L_x$
(so $ L_x$ is a matrix with entries from the localized polynomial ring $\ZZ[x,1/x]$, and $\Theta(x)$ is an integer polynomial in the variables $x$ and $1/x$).
Then
$$
\Theta (x) = K (x-1)^2 x^{-1}  \quad\text{for some integer }K.
$$
Since $\Theta (x) = K (x-1)^2 x^{-1}$, we also have $ K = 2 \Theta(2)$.
\end{Lemma}
 
\noindent Note that the conclusions allow for the possibility that $\Theta$ is identically zero (i.e., $K = 0$).
\begin{proof}
By determinant formulas, it follows that $\det(L_x)$ is a {\it linear} function in the variables $x$ and $x^{-1}$
--- this follows, for example, by looking at the individual terms in the symmetric group sum expansion for a determinant
(see \cite{DF04} Theorem~24, Section~11.4), where only a single factor of $x$ or $1/x$ or $(x)(1/x)$ can appear in each term.
Thus we may write
$$
\Theta(x) = a + b x + c x^{-1} = (x^{-1}) ( a x + b x^2 + c) \qquad \text{for some integers } a,b,c.
$$
Let $q(x)$ be the numerator of the right hand side above. If $q(x) = 0$ then $\Theta(x) = 0$ and so the lemma is true with $K$ chosen to be zero. So assume $q(x)\neq 0.$
Since $\det(L)=0$, it follows that
$x=1$ is a root of $q$, and so all roots of $q$ are rational numbers.
Since $q$ is nonzero but has a root, it is not a constant polynomial.\\
\\
If $s$ is any nonzero root of $q$, then substituting $x = 1/s$ into the (1,2 and 2,1) entries of $L_x$ results in the matrix 
$L_x^{tr}$ (the transpose of $L_x$)
but with $x$ evaluated at $s$ in the latter. Both evaluated matrices have the same determinant, 
hence $1/s$ must be a root of q as well.
Since $q$ has at most two distinct roots, the only possibility for a different root would be $x = -1$. But if $\Theta(x)=K(x^2-1)/x,$ then replacing $x$ by $1/x$ yields $K(1-x^2)x^{-1}$ and so
 $\Theta$ would not satisfy the symmetry condition $\Theta(x) = \Theta(1/x)$.
If $q(x)$ had degree 1, then we would have $\Theta(x) = K (x-1)/x$, for some constant $K$.
But again $\Theta$ would not satisfy the aforementioned symmetry (transpose) condition $\Theta(x) = \Theta(1/x)$,
a contradiction. Thus, $\Theta (x) = K (x-1)^2 x^{-1}$ {for some integer }$K.$

\end{proof}
\begin{Theorem}\label{form}
Assume $X$ satisfies the hypotheses of this subsection and
let $Y$ be the single voltage cover of $X$ by the cyclic group $Z_d = \langle \tau\rangle$.
Then the reduced Stickelberger element may be written in the following form:
$$
\Theta_{Y/X} = K(\tau-1)^2 \tau^{-1}, \qquad\text{for some nonzero integer $K$ independent of $d$.}
$$
\end{Theorem}
\begin{proof}
Let $L_\tau$ be the $n \times n$ voltage Laplacian matrix for the single voltage cover of $X$, so $L_\tau$
has entries in the integral group ring $\ZZ[\tau]$.
The latter ring is the homomorphic image of the polynomial ring $\ZZ[x]$ where the indeterminate $x$ is evaluated at $\tau$.
Since $\tau$ is a unit in the group ring, this algebra homomorphism induces a $\ZZ$-algebra homomorphism
from the localization $\ZZ[x,1/x]$ to $\ZZ[\tau]$ which sends $1/x$ to $\tau^{-1}$ 
(see \cite{DF04}, Theorem~36 and Examples~1 and 2 in Section~15.4).
This further induces a $\ZZ$-algebra homomorphism from the $n \times n$ matrix ring over $\ZZ[x,1/x]$ to
the $n \times n$ matrix ring over $\ZZ[\tau]$.
Since the determinant function is a polynomial in the entries of a matrix, determinants over the former ring map to determinants over the group ring by evaluating $x$ at $\tau$ and $1/x$ at $\tau^{-1}$.
By invoking Lemma~\ref{DetL} and then evaluating $x$ at $\tau$ we obtain that 
\begin{equation}\label{Stickelberger-possobilities}
\Theta_{Y/X} = K (\tau-1)^2 \tau^{-1} , \quad\text{for some integer }K.
\end{equation}

\noindent The entries of $L_x$ depend only on the Laplacian for $X$.
By definition of single voltage cyclic cover, the determinant of $L_\tau$ likewise 
depends only on the Laplacian for $X$ and the choice of generator $\tau$ of the covering group.
The above argument therefore shows that the reduced Stickelberger element is either zero or 
it may be written in a ``canonical form'',
where both the form and the integer constant $K$ are independent of the order, $d$, of the covering group
(although $d$ is inherent in the $\tau$ of this formula).\\
\\
If $\Theta_{Y/X}=0$, then it would follow that $|\mathcal{J}(Y)| = 0$ for all nontrivial (connected) single voltage covering graphs $Y$ of $X$ by Corollary \ref{OrdJac1}, a contradiction.
\end{proof}
%\noindent  \textit{Remark:}\\
%\noindent The reduced Stickelberger element does depend on the choice of generator.  Suppose instead of using generator $\tau$, we use $\mu=\tau^a,$ where $(a,d)=1.$ Then writing the reduced Stickelberger element in terms of $\tau$, we get a different reduced Stickelberger element--namely we get the reduced Stickelberger element 
%$$
%K(\mu-1)^2\mu^{-1}=K(\tau^a-1)^2\tau^{-a}.
%$$
%Note, however, that $K$ does not depend on the choice of generator since for any single voltage cyclic cover, we can write $\Theta_{Y/X}$ uniquely as a $\ZZ$-linear combination of the $\ZZ[G]-$basis elements $1,\tau, \cdots, \tau^{d-1}.$ Then $K$ is found to be the $gcd$ of all coefficients of the linear combination (since if you replace $\tau$ by $\tau^a,$ the $gcd$ of all the coefficients remains the same).
Now let $X$ satisfy the hypothesis at the outset of this subsection and assume $\alpha$ is the single voltage assignment by $Z_d.$ Then by Corollary \ref{OrdJac1}, we have that
\begin{align*}
|\mathcal{J}(Y)|&=\frac{|\mathcal{J}(X)|\prod_{t=1}^{d-1}\frac{|K|(e^{2\pi i t/d}-1)^2}{e^{2\pi i t/d}}}{d}\\
&=\frac{|\mathcal{J}(X)|\cdot |K|^{d-1}d^2}{d}\\
&=|\mathcal{J}(X)|\cdot |K|^{d-1}d,
\end{align*}
\noindent where by the factorization of $\frac{z^d-1}{z-1}$ over $\CC$, the product of all $(e^{2\pi it/d}-1)$-terms in the numerator simplifies to $d^2$ and the product of $e^{2\pi i t/d}$  in the denominator simplifies to~$\pm 1$. 
This proves Corollary \ref{OrdJac2}.

\begin{Corollary}
\hfill
\renum
\item Let $X=K_{n,n}$ with $n\geq 2$ and with the single voltage assignment by $Z_d.$ Then
$$
|\mathcal{J}(Y)|=n^{(2n-4)d+2}(n-1)^{2d-2}d.
$$

\begin{proof}
By Theorem \ref{Kn,n}, we know that the reduced Stickelberger element for these covers is $\Theta_{Y/X}=-(n - 1)^2n^{2n-4}  (\tau - 1)^2  \tau^{-1}.$ 
Since $|\mathcal{J}(X)|=n^{2n-2}$, we get the desired result. 
\end{proof}

\item
Let $X=K_{n,2}$ with $n\geq 2$ and with the directed edge from $v_{1}$ to $v_{n+1}$ with $\tau$ (and so the directed edge from $v_{n+1}$ to $v_{1}$ must be labeled with~$\tau^{-1}$). Then
$$
|\mathcal{J}(Y)|=2^{(n-2)d+1}(n-1)^{d-1}nd.
$$

\begin{proof}
By \cite{phdthesis}, we have that the reduced Stickelberger element is $\Theta_{Y/X}=-2^{n-2}(n-1)(\tau-1)^2\tau^{-1}.$ Since $|\mathcal{J}(X)|=2^{n-1}n,$ we get the desired result.  \end{proof}
\item
Let $X$ be the Petersen graph  with single voltage assignment by $Z_d.$ Then
$$
|\mathcal{J}(Y)|=2^{5d-1}\cdot 5^{2d+1}\cdot d.
$$

\begin{proof}
By \cite{phdthesis}, we know that the reduced Stickelberger element is $\Theta_{Y/X}=-800 (\tau - 1)^2  \tau^{-1}=-2^5\cdot 5^2(\tau-1)^2\tau^{-1}.$ 
Since $|\mathcal{J}(X)|=2^4\cdot 5^3$, we get the desired~result. \end{proof}
\end{enumerate}
\end{Corollary}
%%%%%%%%%%%%%%%%%%%%%%%%%%%%%%%%%%%%%%%%%%%%%%%%%
%\subsection*{Acknowledgements}
%I thank my advisors, Dr. Richard Foote and Dr. Jonathan Sands, for all of their help and inspiration.

%BIBLIOGRAPHY
% You do not have to use the same format for your references, but 
%    include everything in this file.  Don't use natbib please.
% If you use BibTeX to create a bibliography, copy the .bbl file into here.
% \newblock is optional (it adds a little space)

\end{document}